\documentclass[11pt, reqno]{amsart}

\usepackage[margin=3cm]{geometry}
\usepackage{fix-cm}
\usepackage{amsmath, amssymb, amsfonts, amsthm}
\usepackage{mathtools}
\usepackage{graphicx}
\usepackage{algorithm}
\usepackage{algpseudocode}
\usepackage{amsopn}
\usepackage{bm}
\usepackage{multirow}
\usepackage{xcolor}

\usepackage{natbib}

\usepackage{adjustbox}
\usepackage{placeins}
\usepackage{hyperref}

\hypersetup{
  colorlinks = true,
  linkcolor  = blue,
  citecolor  = blue,
  urlcolor   = blue
}

\newtheorem{theorem}{Theorem}
\newtheorem{proposition}{Proposition}
\newtheorem{lemma}{Lemma}
\newtheorem{corollary}{Corollary}
\newtheorem{remark}{Remark}

\DeclareMathOperator{\diag}{diag}
\DeclareMathOperator{\rank}{Rank}
\DeclareMathOperator{\Blockdiag}{Blockdiag}
\DeclareMathOperator{\tr}{Tr}
\newcommand{\R}{\mathbb R}
\newcommand{\N}{\mathbb N}
\newcommand{\x}{\bm x}
\newcommand{\z}{\bm z}
\newcommand{\y}{\bm y}
\newcommand{\e}{\bm e}
\newcommand{\bb}{\bm b}
\newcommand{\ab}{\bm a}

\newcommand{\uu}{\bm u}
\newcommand{\Ss}{\mathcal S}
\newcommand{\PP}{\mathcal P}
\newcommand{\LL}{\mathcal L}

\let\origmaketitle\maketitle
\def\maketitle{
	\begingroup
	\def\uppercasenonmath##1{} 
	\let\MakeUppercase\relax 
	\origmaketitle
	\endgroup
}

\title[Projector-Based Conic Formulations for Distance-Based Assignment]{\Large Projector-Based Conic Formulations and Exact Decompositions for Distance-Based Assignment}

\author[V. Blanco, D. Laborda-Guti\'errez \MakeLowercase{and} M. Mart\'inez-Ant\'on]{
{\large V\'ictor Blanco$^{\dagger}$, Diego Laborda-Guti\'errez$^{\ddagger}$ and Miguel Mart\'inez-Ant\'on$^{\S}$}\medskip\\
$^\dagger$Institute of Mathematics (IMAG), Universidad de Granada\\
\texttt{vblanco@ugr.es}\\
$^\ddagger$Institute of Mathematics (IMUS), Universidad de Sevilla\\
\texttt{dlaborda@us.es}\\
$^\S$Dpt. Mathematics, Statistics \& Computation, Universidad de Cantabria\\
\texttt{mantonm@unican.es}
}


\begin{document}

\begin{abstract}
Given a set of fixed origin points and a set of variable target points in the same ambient metric space, we study distance-based assignment problems in which target placement and the assignment of origins to targets are determined simultaneously to minimize the overall distance-induced cost. We show that the assignment structure can be represented using orthogonal projection matrices. This observation leads to new mixed-projection conic formulations that avoid the big-$M$ envelopes used in standard mixed-integer models. Building on these formulations, we develop exact decomposition algorithms that separate the discrete selection of assignment projectors from the convex conic subproblems for target location. We also derive semidefinite relaxation bounds and primal procedures that provide warm starts for the exact methods. The proposed approach establishes new connections between distance-based assignment, orthogonal projectors, and low-rank semidefinite optimization, and applies to general conic-representable metrics.
\end{abstract}

\maketitle

\section{Introduction}
Assignment problems are ubiquitous in OR/MS, arising in numerous applications such as allocating users to service centers, matching agents to tasks, and distributing resources among projects. In its most general form, an assignment problem is defined as follows. Given a set of \emph{origins} and a set of \emph{targets}. Any origin can be assigned to any target, incurring some cost that may vary depending on the assignment. The objective is to assign each origin to exactly one target so as to minimize the total assignment cost. 

Origins may represent customers, requests, labor, vehicles, products, or data, whereas targets may represent facilities, supply centers, servers, tasks, routes, warehouses, or clusters. Provided that the origins and targets lie in the same ambient metric space, the cost criterion can be based on the distances between them. This is the case, for instance, in facility location and distance-based clustering problems, where decisions on target placement and origin allocation directly affect operating costs, service quality, or the results of data analysis. These problems combine geometric, algebraic, combinatorial, and economic aspects, giving rise to a broad family of models and algorithms.

When the distances between origins and targets are known, the cost function can be represented by a real-valued matrix, whose entries give the assignment costs for each origin-target pair. In this setting, distance-induced costs provide information for decision-making, but the choice of metric space does not affect the methodology. The problem can therefore be viewed as a purely combinatorial optimization problem. This is the case, for instance, in the seminal work on graph medians by \cite{Hakimi1964,Hakimi1965}, who introduced the $p$-median structure for communication networks. In contrast, {we consider problems in which the distance between an origin and a target depends on the assignment itself. Specifically, the problems under study involve a finite set of fixed origin points and a finite set of variable target points, whose locations are determined jointly with the assignment. Throughout this paper, we refer to these problems as \emph{distance-based assignment problems}. Several longstanding problems fall within this framework, including the classical Fermat--Weber problem in industrial location \citep{Weber1909}. In its modern form, the problem seeks the location of a target that minimizes a weighted sum of its distances to the origins.

These problems do not require an artificial set of candidate locations to be specified in advance to construct a matrix of distance-based costs. Instead, the ambient space enters directly into the optimization problem. However, jointly determining the locations of multiple targets and assigning origins to them makes distance-based assignment problems nonconvex and computationally challenging in general.

\subsection{Prior Work}
Combining location and allocation decisions is the main challenge in these problems and has shaped much of the research on the topic. For a single target, optimality conditions and iterative methods exploit convexity and the geometry of sums of distances
\citep{BeckTeboulleChikishev2008,FeketeMitchellBeurer2005}. With multiple targets, the unknown allocation introduces a combinatorial dimension. \cite{Cooper1963} made this structure explicit, and \cite{KuenneSoland1972} later proposed exact and approximate procedures for the multi-target Weber problem. This line of research led to exact enumeration and set-covering methods \citep{Rosing1992}, difference-of-convex global optimization
\citep{ChenHansenJaumardTuy1998}, and geometric branch-and-bound schemes, such
as the big-triangle--small-triangle method
\cite[see, e.g.,][]{BlanqueroCarrizosa2009,DreznerSuzuki2004}.

This difficulty has also motivated growing interest in heuristic methods. \cite{Cooper1964}'s alternating location–allocation algorithm repeatedly assigns each origin to its closest target and then solves an independent single-target problem for each cluster. Projection, neighborhood-search, decomposition, and multistart methods extend this basic idea, enabling larger instances to be tackled \cite[see, e.g.,][]{BrimbergHansenMladenovicTaillard2000}. Exact and heuristic methods play complementary roles. Fast heuristics provide feasible warm starts and starting configurations, while exact methods provide global bounds that help certify solution quality. Recent work has also examined how to construct suitable starting configurations \citep{DreznerBrimberg2026}.

The range of ambient spaces considered in these problems has expanded well beyond the classical planar Euclidean setting. Current models accommodate general $\ell_\tau$ or
heterogeneous metrics, and ordered multitarget criteria
\citep{BlancoPuertoElHaj2014,
BlancoPuertoElHaj2016,BlancoPuertoPonce2017,BlancoGazquezPoncePuerto2023}. This breadth calls for solution methods that are not tied to a particular distance, dimension, or objective. Much of the exact and heuristic literature exploits geometric structure through Voronoi assignments, separability under a fixed allocation, and planar subdivision bounds
\citep{BlanqueroCarrizosa2009,DreznerSuzuki2004,rodriguez2010solving}. Although effective in the planar Euclidean setting, these methods’ bounds and convergence arguments often rely on specific bisectors, subdivisions, or distance properties, hindering their extension to higher dimensions and to general conic-representable metrics. Besides, geometric structure does not eliminate the underlying combinatorial difficulty
\citep{MegiddoSupowit1984}.

Mathematical optimization offers another way to represent this structure. Standard compact formulations use binary variables to encode the allocation and big-$M$ constraints to couple those variables with the distance-based costs. Although these models accommodate broad classes of conic-representable metrics, their convex relaxation bounds could be weak, and selecting valid big-$M$ constants could be difficult. Perspective reformulations can strengthen
mixed-integer nonlinear models with indicator variables
\citep{Gunluk2010-br}, although they still retain an explicit binary
description of the underlying decisions. Semidefinite programming (SDP) provides
a different approach. In combinatorial optimization, matrix liftings have been used to construct convex relaxations of quadratic structures
\citep{Goemans1997semidefinite,VandenbergheBoyd1996}.
\cite{Lasserre2001}'s moment hierarchy extended this idea to polynomial
optimization. Semidefinite relaxations of quadratic problems and low-rank factorization methods have also made this approach useful for a broad range of nonconvex models
\citep{Burer2003-vg,BurerMonteiro2005,LuoMaSoYeZhang2010,Lemon2016}.

The connection is particularly direct in clustering. \cite{peng2005new} introduced a projection-matrix framework for $k$-means-type clustering, and \cite{peng2007approximating} later developed its semidefinite approximation. Subsequent work established recovery guarantees and exactness properties for related conic formulations \citep{li2020birds}.

Rank remains the main source of nonconvexity. Convex surrogates do not, in general, enforce the exact rank associated with a discrete partition, whereas low-rank factorizations reduce the problem dimension but remain nonconvex \citep{Burer2003-vg,BurerMonteiro2005}. Exact methods instead combine convex relaxations with branching or cutting planes. Advances in conic outer approximation have made this approach increasingly practical \citep{BertsimasCoryWright2020,CoeyLubinVielma2020}.

Mixed-projection conic optimization (MPCO), introduced by \cite{BertsimasCoryWrightPauphilet2022}, combines several of these ideas. It represents rank through an
orthogonal projector satisfying $P^2=P$.  This condition is the matrix version
of the binary identity $z^2=z$.  Regularization and conic duality then
separate the choice of a projector from a convex conic subproblem.  This
separation gives globally valid cuts, semidefinite relaxation bounds, and primal
recovery in the same approach.

\subsection{Contribution and Organization}
The purpose of this paper is to study distance-based assignment problems through the lens of orthogonal projectors. An allocation matrix naturally defines an orthogonal projector, allowing the partition structure to be represented directly rather than encoded through products of binary assignment variables and distance-based costs. We show that this observation yields big-$M$-free conic formulations for general conic-representable metrics and leads to exact decomposition algorithms.

The projector-based formulations offer a new perspective on the interplay among assignment, location, rank, and semidefinite optimization. They also provide a foundation for developing formulations and algorithms that do not depend on planar Euclidean geometry.

We derive exact, big-$M$-free conic formulations that couple allocation and location through orthogonal projectors. The sum-aggregated formulation represents an unlabeled partition with a single rank-$p$ projector and enforces a common location within each cluster through a rank constraint.

The decomposition separates partition search from independent convex subproblems for target location. Certified cluster bounds can be reused across
partitions, and semidefinite relaxations supply lower bounds and primal
initializations. We compare the resulting methods with a strengthened
assignment-based conic formulation.

The remainder of the paper is organized as follows. Section \ref{sec:background} introduces the notation and background. Section \ref{sec:pairwise} derives the pairwise projection formulation and its exact logic-based decomposition. Section \ref{sec:mpco} develops the low-rank and mixed-projection formulations and presents an exact outer-approximation algorithm. Section \ref{sec:strengthening} describes the strengthening and heuristic procedures and explains how they improve and initialize the exact solution process. Section \ref{sec:experiments} reports the numerical experiments, and Section \ref{sec:conclusions} concludes the paper.

\section{Distance-Based Assignment and Orthogonal Projectors}\label{sec:background}
Lowercase characters denote scalars, bold lowercase characters denote
vectors, and uppercase characters denote matrices.  For $n\in\N$, we write
$[n]:=\{1,\ldots,n\}$.  The vectors of all ones in $\R^n$ and the identity
matrix in $\R^{n\times n}$ are denoted by $\e_n$ and $I_n$, respectively. Notation $\bm \epsilon_i$ stands for the $i$th vector of the canonical basis.

The origins and targets belong to $\R^m$, which is endowed
with a metric $d$.  We assume that this metric is conic representable.
This means that, for every fixed $\ab\in\R^m$, the epigraph of
$d(\ab,\cdot)$ is the projection of a conic system.  Formally,
there exist an integer $q\geq 0$, a proper cone
$\mathcal{K}$, and an affine map $\LL$ of suitable dimension such that
$$
d(\ab,\x)\leq w
\quad\Longleftrightarrow\quad
\exists\uu\in\R^q:\ \LL(\x,w,\uu)\in\mathcal{K}.
$$
Here, a proper cone is closed, convex, solid, and pointed.  The cone, the
affine map, and the auxiliary dimension may depend on the distance under
consideration.  Every norm-based distance is conic representable.

The class also includes several distances beyond the usual Euclidean and
$\ell_p$ cases.  For example, an anisotropic or Mahalanobis distance can be
written as
$
 d(\ab,\x)=\sqrt{(\x-\ab)^TQ(\x-\ab)}
           =\|L(\x-\ab)\|_2,$
 where $Q=L^TL\succ0,
$
and therefore has a second-order cone representation.  Polyhedral distances
of the form
$
 d(\ab,\x)=\max_{1\leq h\leq H}\{\mathbf{c}_h^T(\x-\ab)\}
$
where the vectors $\mathbf{c}_h$ define a symmetric and point-separating
polytope, are linearly representable by imposing
$\mathbf{c}_h^T(\x-\ab)\leq w$ for every $h\in[H]$.  Distances under refraction,
where travel crosses regions with different transportation costs, are also
conic representable when each regional distance is conic representable.  In
this case, the crossing points are auxiliary variables and the total distance
is the minimum sum of the regional distances
\citep{BlancoPuertoPonce2017}.  Another example is the finite-support Wasserstein distance between two discrete distributions.  For a fixed
ground-cost matrix $C$, it is the optimal value of the transportation problem
$
 \min_{\Pi\geq0}\{\langle C,\Pi\rangle : \Pi\e=\ab,\; \Pi^T\e=\x\},
$
so its epigraph has a linear extended formulation
\citep{cuturi2014fast}.  These examples show that
conic representability may follow from a lifting even when the metric is
not given directly by a standard norm.

For $X\in\R^{m\times n}$, its $j$th column is denoted by
$\x_j\in\R^m$, and its $(i,j)$ entry is denoted by $x_{ij}$.  The
superscript $T$ denotes transposition, and $\tr(X)$ denotes the trace of a
square matrix $X$.  Given matrices $X_i\in\R^{m_i\times n_i}$ for
$i\in[n]$, the matrix $\Blockdiag(X_i)_{i\in[n]}$ has the matrices $X_i$
as its diagonal blocks and zeros elsewhere.  Its dimension is
$\left(\sum_{i=1}^n m_i\right)\times\left(\sum_{i=1}^n n_i\right)$.  The
Euclidean and Frobenius inner products are, respectively,
$
    \langle \x,\y\rangle:=\x^T\y$ and $
    \langle X,Y\rangle:=\tr(X^TY).$

We consider a distance-based assignment problem with origins $\{\ab_1,\ldots,\ab_n\}\subset\R^m$, nonnegative weights
$\omega_i$, and $p$ targets with locations
$\x_1,\ldots,\x_p\in\R^m$.  Each origin is assigned to a closest target.  Since the formulation below requires every target to be assigned to at
least one origin, we assume throughout that
$2\leq p\leq n$.  The distance-based assignment problem is
\begin{equation*}
\min_{X\in\R^{m\times p}} 
\sum_{i=1}^n \omega_i\min_{1\leq j\leq p} d(\ab_i,\x_j).
\end{equation*}
The binary variable $z_{ij}$ takes value one when origin $i$ is
assigned to target $j$.  This gives the mixed-integer nonlinear
formulation:
\begin{align}\tag{MICO}
\min_{X\in\R^{m\times p},\ Z\in \{0,1\}^{n\times p}} &
\sum_{i=1}^{n}\sum_{j=1}^{p} \omega_i z_{ij}d(\ab_i,\x_j) \label{eq:weber_original} \\
\mbox{s.t. }\quad & \sum_{j=1}^{p}z_{ij}=1, && \forall i\in[n]; \label{ctr:1}\\
& \sum_{i=1}^{n}z_{ij}\ge 1, && \forall j\in[p].\label{ctr:2}
\end{align}
When $d$ is induced by an $\ell_p$-norm, the products
$z_{ij}d(\ab_i,\x_j)$ in the objective are usually convexified by combining
$p$-order cone constraints with big-$M$ envelopes
\citep{BlancoGazquezPoncePuerto2023,BlancoPuertoElHaj2016}.  For rational
values of the norm parameter $p$, these $p$-order cones admit extended
representations using second-order cones
\citep[see, e.g.,][]{blanco2025complexity,blanco2024minimal}.  Thus, the resulting models
can be written as mixed-integer second-order cone programs and solved with
available conic optimization solvers.  In this formulation, one
introduces an auxiliary continuous variable for each distance and then
uses the McCormick envelope for the product of this variable and the binary
variable $z_{ij}$.  In an epigraph formulation, this gives constraints of the
form
$
 t_{ij}\geq d(\ab_i,\x_j)-M_i(1-z_{ij})$ and $t_{ij}\geq0.$
The objective becomes $\sum_{i=1}^n\sum_{j=1}^p\omega_i t_{ij}$.  The
constants can be selected separately for each origin as
$
 M_i:=\max_{1\leq k\leq n}d(\ab_i,\ab_k).
$
Indeed, for norm-based distances, a target assigned to a nonempty set of origins can be located in their convex hull without increasing its assigning
cost.  Hence, target location at optimality belongs to the convex hull of
the origins.  By convexity of $d(\ab_i,\cdot)$, its maximum over this
convex hull is attained at a origin, which proves that $M_i$ is a valid
upper bound.

Constraint \eqref{ctr:1} assigns every origin to exactly one target,
and constraint \eqref{ctr:2} ensures that every target is assigned to a nonempty
set of origins.

We next introduce the notation used to reformulate these assignment
decisions.  We let $\Ss^n$ denote the vector space of real $n\times n$
symmetric matrices ($P^T=P$).  A matrix $P\in\R^{n\times n}$ is a projection matrix if
$P^2=P$.  If, in addition, $P\in\Ss^n$, then $P$ is an orthogonal projection
matrix.  The eigenvalues of an orthogonal projector satisfy
$\lambda_i^2=\lambda_i$ and therefore belong to $\{0,1\}$.  Consequently,
$\rank(P)=\tr(P)$.  We use
$
 \PP^n:=\{P\in\Ss^n:P^2=P\}
$
for the set of $n\times n$ orthogonal projection matrices, and
$
 \PP^n_k:=\{P\in\Ss^n:P^2=P,\ \tr(P) =  k\}
$
for the subset of projectors with rank $k$.

Orthogonal projectors provide a matrix representation of linear subspaces.
If $P$ has rank $r$, then it can be written as
$
 P=UU^T,$ with $U^TU=I_r,
$
where the columns of $U\in\R^{n\times r}$ form an orthonormal basis of the
range of $P$.  In particular, every rank-one orthogonal projector has the
form $P=\uu\uu^T$ for some $\uu\in\R^n$ with $\uu^T\uu=1$.  Low-rank
factorizations of this type are widely used to replace a large positive semidefinite
matrix variable by a smaller, but nonconvex, factorized model
\citep[see][]{Burer2003-vg,BurerMonteiro2005,Lemon2016}.

The constraints $P^2=P$ and $\rank(P) = k$ are nonconvex.  A standard
semidefinite relaxation replaces idempotency by the spectral constraints
$0\preceq P\preceq I_n$ together with $\tr(P) = k$.  It is well-known that these constraints
describe the convex hull of orthogonal projectors with rank $k$.
Such semidefinite relaxations have been used to represent partitions in
$k$-means-type models \citep{peng2005new,peng2007approximating}.  In contrast,
mixed-projection conic optimization keeps the projector constraint explicit
and treats it as the matrix counterpart of a binary requirement
\citep{BertsimasCoryWrightPauphilet2022}.  General methods for rank-constrained
and low-rank semidefinite optimization are discussed by
\cite{Lemon2016,SUN2017128}.

Let $Z=(z_{ij})_{i,j}\in\{0,1\}^{n\times p}$ be the assignment matrix, and
let $\z_j$ denote its $j$th column.  Since the rows of $Z$ contain exactly
one nonzero entry, different columns have disjoint supports.  Thus,
$$
Z^TZ
=\diag\left(\z_j^T\z_j\right)_{j\in[p]}
=\diag\left(\e_n^T\z_j\right)_{j\in[p]}.
$$
All diagonal entries are positive by \eqref{ctr:2}, so $Z^TZ$ is invertible.
We can therefore define
$
P:=
Z(Z^T Z)^{-1}Z^T \in \PP^n\cap \R_+^{n\times n},
$
which is the orthogonal projector onto the column space of $Z$.  This matrix
is entrywise nonnegative.  Constraint \eqref{ctr:1} can be written as
$Z\e_p=\e_n$.  Hence, $\e_n$ belongs to the column space of $Z$, so that
$P\e_n=\e_n$.  Since all $p$ columns of $Z$ are nonzero and mutually
orthogonal, $\rank(Z)=p$.  Therefore, $\tr(P)=\rank(P)=p$.

The projector $P$ can also be disaggregated by target as
$$
P=\sum_{j=1}^pP_j,
\qquad
P_j:=
\frac{\z_j \z_j^T}{\e_n^T\z_j}.
$$
For every $j\in[p]$, $P_j\in\PP^n_1\cap\R_+^{n\times n}$ is the rank-one
orthogonal projector associated with the origins assigned to target
$j$.  In particular, $P_j\e_n=\z_j$.  These identities will be used in the
projector-based formulations developed in the subsequent sections.

\section{Pairwise Formulation}
\label{sec:pairwise}

We first derive a projector-based conic formulation directly from the assignment model
\eqref{eq:weber_original}.  The main idea is to replace the indicator vector
of each target by the rank-one projector onto the subspace that it spans.
To recover the contribution of a cluster from this normalized projector, we
write its service cost as an average of pairwise quantities.  This produces
an exact formulation in which the assignment decisions are carried only by
nonnegative orthogonal projectors.  Unlike a standard mixed-integer conic
formulation, the resulting model does not require big-$M$ constants.

For every $\x\in\R^m$, define the symmetric matrix
$Y(\x)=(y_{ik}(\x))_{i,k\in[n]}\in\Ss^n$ componentwise by
$
 y_{ik}(\x):=\omega_i d(\ab_i,\x)+\omega_k d(\ab_k,\x).
$
The entry $y_{ik}(\x)$ is the total weighted cost of assigning origins
$i$ and $k$ to a same target located at $\x$.  Hence, $Y(\x_j)$
collects all pairwise assignment costs associated with target $j$. We obtain
\begin{align}
 \sum_{i=1}^{n}\omega_i z_{ij}d(\ab_i,\x_j)
 &=\frac{1}{2\e_n^T\z_j}
   \sum_{i=1}^{n}\sum_{k=1}^{n}
   z_{ij}z_{kj}y_{ik}(\x_j) \nonumber\\
 &=\frac12\langle P_j,Y(\x_j)\rangle.
 \label{eq:pairwise-identity}
\end{align}
The factor $1/2$ corrects the two occurrences of every individual assignment
cost in $y_{ik}(\x_j)$.  Identity \eqref{eq:pairwise-identity} is the link
between the classical assignment objective and the projector-based model.

\begin{theorem}\label{th:pw}
The distance-based assignment problem has the following exact formulation:
\begin{align}
 \min\quad &
 \frac12\sum_{j=1}^{p}\langle P_j,Y(\x_j)\rangle
 \tag{PW}\label{mod:pairwise}\\
 \textnormal{s.t.}\quad
 &\sum_{j=1}^{p}P_j\e_n=\e_n,\nonumber\\
 &P_j\in \PP^n_1\cap\R_+^{n\times n},\quad \x_j\in\R^m, \quad \forall j\in[p].\nonumber
\end{align}
The assignment vectors of any feasible solution of \eqref{mod:pairwise} are
recovered as $\z_j=P_j\e_n$, $j\in[p]$.  Conversely, every feasible
assignment matrix $Z$ induces feasible projectors through
\begin{equation}
 P_j=\frac{\z_j\z_j^T}{\e_n^T\z_j},\qquad j\in[p].
 \label{eq:assignment-to-pairwise-projector}
\end{equation}
Both transformations preserve the objective value.
\end{theorem}

\begin{proof}
Let $Z$ be a feasible assignment matrix.  Constraint \eqref{ctr:2} implies
$\e_n^T\z_j>0$ for every $j\in[p]$, so the matrices in
\eqref{eq:assignment-to-pairwise-projector} are well defined.  Each $P_j$ is
entrywise nonnegative, symmetric, idempotent, and has rank one.  Also,
$P_j\e_n=\z_j$.  It follows from \eqref{ctr:1} that
$\sum_{j=1}^pP_j\e_n=\e_n$, and
\eqref{eq:pairwise-identity} proves that the objective value is preserved.

For the reverse implication, consider a feasible collection
$(P_j)_{j\in[p]}$ in \eqref{mod:pairwise} and set
$\z_j=P_j\e_n$.  Nonnegativity and the coupling constraint give
$0\leq z_{ij}\leq1$ and $\sum_{j=1}^p z_{ij}=1$ for every $i\in[n]$.
Idempotency and symmetry yield
$
 \z_j^T\z_j
 =\e_n^TP_j^TP_j\e_n
 =\e_n^TP_j\e_n
 =\e_n^T\z_j.
$
Consequently,
$$
 0=\sum_{j=1}^p\bigl(\e_n^T\z_j-\z_j^T\z_j\bigr)
  =\sum_{i=1}^n\sum_{j=1}^p z_{ij}(1-z_{ij}).
$$
Every term in the last sum is nonnegative.  Therefore,
$z_{ij}\in\{0,1\}$ for all $i\in[n]$ and $j\in[p]$.  The linear equality excludes the zero projector, so each induced cluster is
nonempty.  Finally, a nonnegative rank-one projector whose row-sum vector is
$\z_j$ must satisfy
$P_j=\z_j\z_j^T/(\e_n^T\z_j)$.  Applying
\eqref{eq:pairwise-identity} completes the proof.
\end{proof}

We next express \eqref{mod:pairwise} in the value-function form used in
mixed-projection conic optimization.  Define the epigraph set
\begin{equation}
 \mathcal E:=\left\{(W,\x)\in\Ss^n\times\R^m:
 w_{ik}\geq y_{ik}(\x),\ i,k\in[n]\right\}
 \label{eq:pairwise-epigraph-set}
\end{equation}
and let $\delta_{\mathcal E}$ be its extended-valued indicator function.
Thus, $\delta_{\mathcal E}(W,\x)=0$ when
$(W,\x)\in\mathcal E$ and $+\infty$ otherwise. Since the metric $d$ is
conic representable, set \eqref{eq:pairwise-epigraph-set} is also conic
representable because each of its inequalities is an epigraph constraint for
a nonnegative weighted sum of two distances.

For $P\in\PP_1^n\cap\R_+^{n\times n}$, consider the extended-valued
function
\begin{equation}
 f(P):=\inf_{V,W,\x}
 \left\{\tr(V)+\delta_{\mathcal E}(W,\x):
 V=PW\right\}.
 \label{eq:pairwise-value-function}
\end{equation}
The matrix $V$ isolates the nonconvexity stemming from the product between the projector and the epigraph variable.  For fixed $P$, the equality $V=PW$ is linear, so
\eqref{eq:pairwise-value-function} is a conic optimization problem.

\begin{theorem}
For every $P\in\PP_1^n\cap\R_+^{n\times n}$,
\begin{equation}
 f(P)=\inf_{\x\in\R^m}\langle P,Y(\x)\rangle.
 \label{eq:pairwise-value-identity}
\end{equation}
Consequently, \eqref{mod:pairwise} is equivalent to
\begin{align}
 \min\quad &
 \frac12\sum_{j=1}^p f(P_j)
 \tag{PW-VF}\label{mod:pairwise-value-function}\\
 \textnormal{s.t.}\quad
 &\sum_{j=1}^pP_j\e_n=\e_n,\nonumber\\
 &P_j\in\PP_1^n\cap\R_+^{n\times n}, \quad \forall j\in[p].\nonumber
\end{align}
In particular, \eqref{mod:pairwise-value-function} separates the discrete
partition structure, represented by the projectors, from the location conic
subproblems represented by $f$.
\end{theorem}

\begin{proof}
Fix $P\in\PP_1^n\cap\R_+^{n\times n}$.  Eliminating $V$ from
\eqref{eq:pairwise-value-function} gives
$$
 f(P)=\inf_{(W,\x)\in\mathcal E}\tr(PW)
     =\inf_{(W,\x)\in\mathcal E}\langle P,W\rangle,
$$
where the latter equality follows because $P$ and $W$ are symmetric.  For any
$(W,\x)\in\mathcal E$, the matrix $W-Y(\x)$ is entrywise nonnegative.
Since $P$ is also entrywise nonnegative,
$
 \langle P,W\rangle\geq\langle P,Y(\x)\rangle.
$
Equality is attained by choosing $W=Y(\x)$.  Taking the infimum over $\x$
proves \eqref{eq:pairwise-value-identity}.  Substitution of this identity in
\eqref{mod:pairwise} proves the equivalence.
\end{proof}

The value-function representation does not by itself convexify the idempotency requirement.  It does, however, separate the selection of a partition from
the continuous location decisions.  This separation leads directly to the
following exact decomposition method and completes our analysis of the
pairwise projector-based model.

\subsection{Exact Logic-Based Decomposition}
Let $\PP^{n,p}$ define as
\begin{equation}
 \left\{(P_1,\ldots,P_p)\in
 (\PP_1^n\cap\R_+^{n\times n})^p :
 \sum_{j=1}^pP_j\e_n=\e_n\right\}.
 \label{eq:pairwise-projector-set}
\end{equation}
By Theorem~\ref{th:pw}, $\mathcal P^{n,p}$ is precisely the
set of labeled nonempty partitions represented through orthogonal
projectors.  For a fixed projector collection
$\mathbf{P}=(P_1,\ldots,P_p)\in\mathcal P^{n,p}$, define the location
subproblem value
\begin{equation}
 g(\mathbf{P}):=
 \min_{X\in \R^{m\times p}}
 \frac12\sum_{j=1}^p\langle P_j,Y(\x_j)\rangle
 =\frac12\sum_{j=1}^p f(P_j).
 \label{eq:pairwise-subproblem-value}
\end{equation}
The equality follows from \eqref{eq:pairwise-value-identity}.  Thus,
$g(\mathbf{P})$ is obtained by solving $p$ independent conic programs.

To construct a master problem, consider an evaluated projector collection
$\bar{\mathbf{P}}=(\bar P_1,\ldots,\bar P_p)$, set
$\bar\z_j=\bar P_j\e_n$, and define
\begin{equation}
 \Delta_{\bar{\mathbf{P}}}(\mathbf{P}):=
 \sum_{j=1}^p\bar\z_j^T(\e_n-P_j\e_n).
 \label{eq:partition-distance}
\end{equation}
For $\mathbf{P}\in\mathcal P^{n,p}$, the vector $P_j\e_n$ is the binary
incidence vector of cluster $j$.  Consequently,
$\Delta_{\bar{\mathbf{P}}}(\mathbf{P})=0$ if $\mathbf{P}$ and
$\bar{\mathbf{P}}$ represent the same labeled partition, whereas
$\Delta_{\bar{\mathbf{P}}}(\mathbf{P})\geq1$ for every different labeled
partition.  Since $g(\mathbf{P})\geq0$, the logic-based optimality cut
\begin{equation}
 \theta\geq
 g(\bar{\mathbf{P}})
 \bigl(1-\Delta_{\bar{\mathbf{P}}}(\mathbf{P})\bigr)
 \label{ctr:pairwise-logic-cut}
\end{equation}
is valid over $\mathcal P^{n,p}$.  It imposes the exact location cost at
$\bar{\mathbf{P}}$ and becomes redundant at every other labeled partition.
After evaluating the projector collections
$\mathbf{P}_1,\ldots,\mathbf{P}_r$, the master problem is
\begin{align}
 \min_{\mathbf{P}\in\mathcal P^{n,p},\theta \geq 0}\quad &\theta
 \tag{PW-LB$_r$}\label{mod:pairwise-benders-master}\\
 \textnormal{s.t.}\quad
 &\theta\geq g(\mathbf{P}_\ell)
  \bigl(1-\Delta_{\mathbf{P}_\ell}(\mathbf{P})\bigr), \quad \forall \ell\in[r].\nonumber
\end{align}

\begin{proposition}
The optimal value of the master problem
\eqref{mod:pairwise-benders-master} is a lower bound on that of
\eqref{mod:pairwise}.  If a master solution not
evaluated previously is selected at every iteration, then the procedure that
alternates between \eqref{mod:pairwise-benders-master} and the location
subproblem \eqref{eq:pairwise-subproblem-value} terminates after finitely many
iterations with a globally optimal solution of \eqref{mod:pairwise}.
\end{proposition}

\begin{proof}
Cut \eqref{ctr:pairwise-logic-cut} never overestimates the value of a
partition.  Indeed, its right-hand side equals
$g(\bar{\mathbf{P}})$ at the evaluated partition and is nonpositive at any
other partition, while all location costs are nonnegative.  Therefore, the
master is a relaxation of the value-function problem and supplies a valid
lower bound.  The set $\mathcal P^{n,p}$ represents finitely many labeled
partitions.  Each iteration evaluates a new member of this set and adds a cut
that reproduces its exact value.  Hence, after finitely many iterations,
either the master lower bound equals the best evaluated value or every
partition has been evaluated.  In both cases, the best evaluated solution is
globally optimal.
\end{proof}

This master can be strengthened by exploiting the
target-wise separability of $g(\mathbf P)$ \eqref{eq:pairwise-subproblem-value}.  It introduces
variables $\theta_j\geq0$, imposes $\theta\geq\frac{1}{2}\sum_{j=1}^p\theta_j$, and minimizes
$\theta$.  Let $C\subset [n]$ represent the cluster of origins induced by a previously evaluated rank-one projector $\bar P_j$. If $L(C)$ denotes a dual-certified value for the single-target problem $f(\bar P_j)$, it adds
\begin{equation}\label{eq:sep-cuts}
    \theta_j\geq L(C)\left(\sum_{i\in C}\bm\epsilon_i^TP_j\e_n-|C|+1\right),
 \qquad j\in[p].
\end{equation}
This inequality remains valid when target $j$ contains $C$, so a certified
cluster value can strengthen many partitions and every target label.  The
aggregate cut \eqref{ctr:pairwise-logic-cut} is retained.  The program can be written as a nonconvex quadratically constrained problem
by factoring $P_j=\uu_j\uu_j^T$.  Any global optimization solver able to
handle this master can be used.  At each master incumbent, the location
subproblem \eqref{eq:pairwise-subproblem-value} can be handled by any suitable
conic solver, and cut \eqref{ctr:pairwise-logic-cut} can be added through a
lazy-constraint callback.  Alternatively, an explicit solve--separate--resolve
loop can solve the master, certify the selected partition through the conic
subproblem, add the resulting cuts, and solve the strengthened master again. 
The best evaluated value is an upper bound, while the master
value is a lower bound.  The method therefore provides a global optimality
certificate when these bounds coincide.  The cuts are intentionally simple
and can be weak.  The warm-start procedures described in Section
\ref{sec:strengthening}, together with symmetry-breaking conditions for target
labels, can substantially reduce the number of partitions evaluated.

    \begin{remark}
    Notice that a similar exact decomposition scheme, developed above to tackle the mixed-projection program \eqref{mod:pairwise}, can be adapted to the mixed-integer program \eqref{eq:weber_original}. We can express \eqref{eq:weber_original} in the value-function form. Define the epigraph set
    $$
        \widetilde{\mathcal E}:= \{(W,X)\in \R^{n\times p}\times \R^{m\times p} : w_{ij}\geq \omega_i d(\ab_i,\x_j),\; i\in [n], j\in [p]\}.
    $$
    For an assignment matrix $Z\in \{0,1\}^{n\times p}$ consider the extended-valued function
    \begin{equation*}
        \widetilde{f}(Z) := \inf_{W,X} \{\tr(V) + \delta_{\tilde{\mathcal E}}(W,X) : V=Z^TW\}.
    \end{equation*}
    Now, the matrix $V$ isolates the product between binary variables and the continuous epigraph variables. For fixed $Z$, the equality $V=Z^TW$ is linear, so $\widetilde{f}(Z)$ is a conic optimization problem which serves as a tractable subproblem in the decomposition scheme. It is not difficult to see that 
    \begin{equation*}
    \min_{Z\in\{0,1\}^{n\times p}} \{\widetilde{f}(Z) : Z\e_p = \e_n, \; Z^T\e_n \ge \e_p\}
    \end{equation*}
    is equivalent to \eqref{eq:weber_original}. Hence, we can add the logic-based optimality cuts
    \begin{equation*}
        \theta\geq \widetilde{f}(\bar{Z})\left(\tr(\bar{Z}^TZ)+1 -\e_n^T\bar{Z}\e_p \right)
    \end{equation*}
    provided by the subproblem evaluated at a fixed assignment $\bar{Z}$ in order to construct the master problem as detailed above. This approach handles nonconvexity directly within the decomposition framework, without requiring the model to be convexified through big-$M$ techniques, as is commonly done in the state-of-the-art literature. However, this remains an oracle-based reformulation of the classical assignment model.  Each cut is informative only at its evaluated assignment, so all labeled assignments may still need to be examined.  The master also retains $np$ binaries and target-label symmetry, and each new assignment requires a conic subproblem.  In contrast, the projector formulations admit semidefinite relaxation bounds, projection inequalities, trace bounds, and matrix cuts that act on several partitions, while the sum of rank-one projectors is label invariant.  Thus, eliminating big-$M$ constants from the assignment model does not provide the same lower-bounding and symmetry-reduction mechanisms. To the best of our knowledge, this approach has not yet been explored in the context of distance-based assignment. Although this idea is not the focus of the present paper, it represents a promising direction for future research.
\end{remark}

The next section develops a lifted formulation that represents the
projector-location interaction through conic and semidefinite constraints.

\section{Low-Rank Formulations}
\label{sec:mpco}

The pairwise model developed in Section \ref{sec:pairwise} gives an exact
projector-based representation of the assignment decisions.  Its objective still
couples each projector with a matrix of distances that depends on an unknown
target location.  We now derive a second representation that moves this
interaction from the objective to a projection constraint.  The derivation
has two steps.  We first use the perspective of each distance to represent
an active assignment-location term without a big-$M$ constant.  We then
relate the assignment and location variables through orthogonal projectors.  The result is
an exact mixed-projection conic optimization model whose continuous part is
conic once the projectors are fixed.

Perspective formulations strengthen mixed-integer convex models with
indicator variables \citep{Gunluk2010-br}.  Low-rank factorizations and
rank-constrained semidefinite models have been studied extensively
\citep{Burer2003-vg,BurerMonteiro2005,Lemon2016,SUN2017128}.  Our construction
combines these ideas within the low-rank MPCO framework of
\cite{BertsimasCoryWrightPauphilet2022}.  In contrast with
the standard assignment formulation, it requires neither binary variables
nor big-$M$ envelopes.

Given a feasible assignment and its centers, define
\begin{equation}
 Q_j:=\z_j\x_j^T\in\R^{n\times m}, \quad \forall j\in [p].
 \label{eq:Q-definition}
\end{equation}
The active distance can be written using only $P_j$ and $Q_j$.  For this
purpose, let $d_i(\x):=d(\ab_i,\x)$ and denote its recession function by
$d_i^\infty$.  We define the closed perspective of $d_i$ as
\begin{equation*}
 \widehat d_i(\alpha,\bm q):=
 \begin{cases}
  \alpha d_i(\bm q/\alpha),&\alpha>0,\\
  d_i^\infty(\bm q),&\alpha=0,\\
  +\infty,&\alpha<0.
 \end{cases}
\end{equation*}
This definition depends only on the arguments $(\alpha,\bm q)$ and does not
involve an assignment variable.  In particular,
$\widehat d_i(1,\x)=d(\ab_i,\x)$ and
$\widehat d_i(0,\bm 0)=d_i^\infty(\bm 0)=0$.  For a feasible assignment,
$P_j\e_n=\z_j$ and $Q_j^T\bm \epsilon_i=z_{ij}\x_j$.  Therefore,
\begin{equation}
 \widehat d_i\bigl(\bm\epsilon_i^TP_j\e_n,Q_j^T\bm \epsilon_i\bigr)
 =z_{ij}d(\ab_i,\x_j).
 \label{eq:distance-identity}
\end{equation}
In formulations \eqref{mod:mpco-disaggregated}, and \eqref{mod:mpco-aggregated}, the two
arguments of $\widehat d_i$ are the decision expressions
$\bm\epsilon_i^TP_j\e_n$ and $Q_j^T\bm \epsilon_i$.  Thus, neither $z_{ij}$ nor $\x_j$ is a
decision variable in these models.  They are used only in
\eqref{eq:distance-identity} to prove equivalence with the original
formulation.  The epigraph of $\widehat d_i$ is obtained by homogenizing an
epigraph representation of $d_i$, and is therefore conic representable
whenever the epigraph of $d_i$ is conic representable
\citep{Gunluk2010-br}.

More specifically, suppose that a conic representation of $d_i$ is
$$
 d_i(\x)\leq w
 \quad\Longleftrightarrow\quad
 \exists\,\uu:\ A_i\x+\bb_iw+C_i\uu+\bm c_i\in\mathcal K_i,
$$
Homogenizing this equivalence gives
$$
 t_{ij}\geq
 \widehat d_i\bigl(\bm\epsilon_i^TP_j\e_n,Q_j^T\bm \epsilon_i\bigr)
 \Longleftrightarrow
 \left\{
 \begin{array}{l}
 \bm\epsilon_i^TP_j\e_n\geq0,\\
 A_iQ_j^T\bm \epsilon_i+\bb_i t_{ij}+C_i\bm v_{ij}
       +\bm\epsilon_i^TP_j\e_n\bm c_i\in\mathcal K_i
 \end{array}
 \right.
$$
for an auxiliary vector $\bm v_{ij}$.  All expressions above are affine in
$(P_j,Q_j,t_{ij},\bm v_{ij})$.  Hence, this is a standard conic constraint and does not introduce
a product among decision variables.  For example, under the Euclidean
distance, it reduces to the second-order cone constraint
$
 \|Q_j^T\bm \epsilon_i-\bm\epsilon_i^TP_j\e_n\ab_i\|_2\leq t_{ij}.
$
For a Mahalanobis distance with factor $L$, the left-hand side becomes 
$\|L(Q_j^T\bm \epsilon_i-\bm\epsilon_i^TP_j\e_n\ab_i)\|_2$.  For a polyhedral distance, the
corresponding implementation consists only of linear inequalities.  Thus,
an implementation creates the scalar expression $\bm\epsilon_i^TP_j\e_n$, the vector
expression $Q_j^T\bm \epsilon_i$, and $t_{ij}$, and submits their joint conic epigraph
directly to the conic solver.  However, this implementation needs to
enforce consistency between the first two expressions, that is, to enforce that $Q_j = \bm z_j \bm x_j^T$.

We next use the central observation of mixed-projection conic
optimization: a rank restriction can be represented by requiring a matrix
to be invariant under a low-rank orthogonal projector
\citep{BertsimasCoryWrightPauphilet2022}.

\begin{lemma}
\label{prop:rank-projector}
For any $W\in\R^{r\times s}$ and $k\in[r]$, $
 \rank(W) \leq k$ if and only if there exists $P\in\PP^r_k$ such that $W=PW$. Furthermore, if $P = \sum_{j=1}^k v_j v_j^T$ for orthonormal vectors
$v_1,\dots,v_k$, then $W=PW$ is equivalent to
$W = \sum_{j=1}^k v_j w_j^T$ for some vectors
$w_1,\dots,w_k \in \mathbb R^s$.
\end{lemma}

Applying this lemma, $Q_j = \bm z_j \bm x_j^T$ can be enforced and the following result can be proved.

\begin{theorem}
The distance-based assignment problem is equivalent to
\begin{align*} 
 \min\quad
 &\sum_{j=1}^p\sum_{i=1}^n\omega_i t_{ij}
 \tag{LR-D}\label{mod:mpco-disaggregated}\\
 \textnormal{s.t.}\quad
 &\sum_{j=1}^pP_j\e_n=\e_n,\nonumber\\
 &t_{ij}\geq
   \widehat d_i\bigl(\bm\epsilon_i^TP_j\e_n,Q_j^T\bm \epsilon_i\bigr),
   &&i\in[n],\ j\in[p];\nonumber\\
 &Q_j=P_jQ_j, &&j\in[p];\nonumber\\
 &P_j\in\PP^n_1 \cap \R^{n \times n}_+, &&j\in[p].\nonumber
\end{align*}
\end{theorem}

\begin{proof}
Given a clustering, construct $P_j$ and $Q_j$, as in
\eqref{eq:Q-definition}.  Each $Q_j$ is such that \eqref{eq:distance-identity} gives a feasible point with the same objective
value while satisfying $Q_j = P_j Q_j$, since
$$ P_j Q_j = \frac{\bm z_j \bm z_j^T}{\bm e_n^T \bm z_j} \bm z_j \bm x_j^T = \frac{\bm z_j^T \bm z_j}{\bm e_n^T \bm z_j} \bm z_j \bm x_j^T = \bm z_j \bm x_j^T = Q_j. $$
Conversely, constraints $P_j \in \PP_1^n \cap \R^{n\times n}_+$, together with $\sum_{j=1}^pP_j\e_n=\e_n$, induce the
incidence vectors $\z_j=P_j\e_n$ of a partition of $[n]$. Using the preceding lemma, constraints
$Q_j = P_j Q_j$ imply $Q_j=\z_j\x_j^T$ for some $\x_j\in\R^m$. Therefore, the conic
inequalities reduce to $t_{ij}\geq z_{ij}d(\ab_i,\x_j)$, and they are tight
at optimality.
\end{proof}

For an aggregated representation, write
$
  \widetilde{P}:=\Blockdiag(P_1,\ldots,P_p) \in \R^{pn \times pn}$ and $\widetilde{Q}:=\Blockdiag(Q_1,\ldots,Q_p) \in\R^{pn\times pm}.$
Because of this block structure, it is clear that these aggregated matrices also satisfy $\widetilde{Q} = \widetilde{P} \widetilde{Q}$ while $\widetilde{P} \in \PP_p^{pn}$. 

\begin{corollary} The distance-based assignment problem is equivalent to the
following formulation:
\begin{align}
 \min\quad
 &\sum_{j=1}^p\sum_{i=1}^n\omega_i t_{ij}
 \tag{LR-A}\label{mod:mpco-aggregated}\\
 \textnormal{s.t.}\quad
 &\sum_{j=1}^pP_j\e_n=\e_n,\nonumber\\
 &\tr(P_j)=1, &&j\in[p],\nonumber\\
 &t_{ij}\geq
   \widehat d_i\bigl(\bm\epsilon_i^TP_j\e_n,Q_j^T\bm \epsilon_i\bigr),
   &&i\in[n],\ j\in[p];\nonumber\\
 &\widetilde{Q}= \widetilde{P} \widetilde{Q};\nonumber\\
 &\widetilde{P} \in\PP^{pn}_p \cap \R^{pn \times pn}_+.\nonumber
\end{align}
Formulation \eqref{mod:mpco-disaggregated} uses one rank-one projector for each target,
whereas this formulation uses one rank-$p$ projector for the
block-diagonal matrix $\widetilde{P}$.
\end{corollary}

\begin{proof}
Equivalence follows from the block structure of $\widetilde{P}$ and $\widetilde{Q}$. We have that $Q_j = P_j Q_j$ and, since $\tr(P_j) = 1$, each block $P_j$ is non-zero and $P_j \in \PP^n_1$.
Then, this formulation is equivalent to \eqref{mod:mpco-disaggregated} and therefore to the distance-based assignment problem.
\end{proof}

A third possible low-rank model arises by aggregating the matrices through a sum instead of a diagonal arrangement. Since the $\bm z_j$ vectors of a valid assignment are orthogonal, all the information about the assignment is still preserved in this rank $p$ sum, resulting in a more compact reformulation.

\begin{lemma}\label{lem:aggregate-partition-projector}
Let $P\in\PP_p^n\cap\R_+^{n\times n}$ satisfy $P\e_n=\e_n$.  Then there is
a partition $\{C_1,\ldots,C_p\}$ of $[n]$ with incidence vectors
$\z_1,\ldots,\z_p\in\{0,1\}^n$ such that $P=\sum_{j=1}^p\frac{\z_j\z_j^T}{\e_n^T\z_j}$. 
\end{lemma}
\begin{proof}
The connected components of the positive-entry graph of the nonnegative
symmetric matrix $P$ permit a simultaneous permutation of its rows and
columns into irreducible diagonal blocks.  Because $P\e_n=\e_n$ and
$P^2=P$, every block is stochastic and idempotent.  Perron--Frobenius theory
implies that eigenvalue one is simple in each block, while idempotency implies
that every eigenvalue is zero or one.  Thus, every block has rank one and,
by symmetry and stochasticity, equals $|C|^{-1}\e_{|C|}\e_{|C|}^T$ on its
index set $C$.  After undoing the permutation, this block is
$\z\z^T/(\e_n^T\z)$, where $\z$ is the incidence vector of $C$.  Finally,
$\rank(P)=p$ implies that there are exactly $p$ blocks.
\end{proof}

\begin{theorem}
    The following formulation is equivalent to the distance-based assignment problem:
    \begin{align}
 \min\quad
 &\sum_{i=1}^n\omega_i t_i
 \tag{LR-SA}\label{mod:mpco-sum-aggregated}\\
 \textnormal{s.t.}\quad
 &t_i\geq
   d\bigl(\ab_i,Q^T\bm \epsilon_i\bigr), \quad i\in[n];\nonumber\\
 & P \bm e_n = \bm e_n, \nonumber\\
 & Q=  P  Q,\nonumber\\
 & P \in\PP^{n}_p \cap \R^{n \times n}_+,  Q \in \mathbb R^{n \times m}.\nonumber
\end{align}
\end{theorem}
\begin{proof}
    Given a clustering, construct the matrices $P_j$ and $Q_j$ as before and choose $P = \sum_{j=1}^p P_j$ and $Q = \sum_{j=1}^p Q_j$. Such a choice has $P \in \PP^n_p \cap \mathbb R^{n \times n}_+$. Since vectors $\bm z_j$ are pairwise orthogonal, we have that
    $$
    PQ = \sum_{j,j'=1}^p P_j Q_{j'} = \sum_{j,j'=1}^p \frac{\bm z_j \bm z_j^T}{\bm e_n^T \bm z_j} \bm z_{j'} \bm x_{j'}^T = \sum_{j=1}^p \frac{\bm z_j^T \bm z_j}{\bm e_n^T \bm z_j} \bm z_j \bm x_j^T = \sum_{j=1}^p Q_j = Q.
    $$
    Furthermore, $Q^T \bm \epsilon_i = \sum_{j=1}^p Q_j^T \bm \epsilon_i = \sum_{j=1}^p \bm x_j \bm z_j^T \bm \epsilon_i = \sum_{j=1}^p z_{ij} \bm x_j$, 
    that is, $Q^T \bm \epsilon_i$ is the center assigned to $\bm a_i$. Therefore, this construction gives a feasible solution with the same objective value, and the inequalities $t_i \geq d(\bm a_i, Q^T \bm \epsilon_i)$ are tight at optimality.

    Conversely, given any feasible matrices $P$ and $Q$ of this formulation,
    Lemma~\ref{lem:aggregate-partition-projector} gives
    $P=\sum_{j=1}^p \bm z_j\bm z_j^T/(\bm e_n^T\bm z_j)$ for the incidence
    vectors of a partition. Lemma~\ref{prop:rank-projector} assures $Q$ is of the form $Q = \sum_{j=1}^p Q_j = \sum_{j=1}^p \bm z_j \bm x_j^T$ for some vectors $\bm x_j \in \mathbb R^m$. Again, $Q^T \bm \epsilon_i = \sum_{j=1}^p z_{ij} \bm x_j$ will act as the center assigned to $\bm a_i$ and constraints $t_i \geq d(\bm a_i, Q^T \bm \epsilon_i)$ will be tight at optimality, proving the equivalence of the problems.
\end{proof}

Once the projectors are fixed, all remaining constraints in these models are
conic. Thus, the assignment decisions of the original big-$M$ model have
been replaced by orthogonal projection matrices.

\begin{remark} Related semidefinite and conic
decomposition techniques are discussed by
\cite{BertsimasCoryWright2020,CoeyLubinVielma2020,VandenbergheBoyd1996}. In fact, the assignment matrices $P$ and the center matrices $Q$ in the above three formulations play the same role as matrices $Y$ and $X$ respectively in the framework described by \cite{BertsimasCoryWright2020}.

Furthermore, it is of note how well projection matrices seem to be suited to the role of assignment variables, particularly in the context of clustering and distance-based assignment problems. The earlier work of \cite{peng2005new} on semidefinite formulations for $k$-means-type problems already hinted at this interplay.
\end{remark}

\subsection{Exact Outer Approximations}
We close the section with a globally exact method for the disaggregated model
\eqref{mod:mpco-disaggregated}.  The method follows the outer-approximation
strategy for MPCO in \cite{BertsimasCoryWrightPauphilet2022}.  A global
projector master selects the partition structure, while conic subproblems
evaluate the location decisions and generate supporting cuts.
The same construction can be applied to the sum-aggregated
formulation \eqref{mod:mpco-sum-aggregated}.  We first present the
disaggregated version and then state its sum-aggregated counterpart.



Let $\mathcal C$ collect the conic constraints in
\eqref{mod:mpco-disaggregated} involving $(P,Q,T)$, excluding
$Q_j=P_jQ_j$.  For
$\mathbf P=(P_1,\ldots,P_p)\in\mathcal P^{n,p}$, define the fixed-projector
value:
$$
 h(\mathbf P):=
 \min_{Q,T}
 \left\{
  \sum_{j=1}^p\sum_{i=1}^n\omega_i t_{ij}:
  (P,Q,T)\in\mathcal C,\ Q_j=P_jQ_j\ j\in[p]
 \right\}.
$$
The subproblem is conic for fixed $\mathbf P$ and can be solved with any
suitable conic optimization solver.  It
contains neither binary variables nor big-$M$ convexifications.

Since
$h(\mathbf P)\geq0$, the projector-space logic cut
\begin{equation}
 \theta\geq h(\bar{\mathbf P})
 \bigl(1-\Delta_{\bar{\mathbf P}}(\mathbf P)\bigr)
 \label{ctr:optimality-cut}
\end{equation}
is valid (recall $\Delta$ as \eqref{eq:partition-distance}).  Because \eqref{eq:pairwise-projector-set} contains only partition
projectors and every such projector has a conic completion, no
feasibility cuts or artificial bounds are required.

After $r$ evaluations, the master problem reads
\begin{align}
 \min_{\mathbf P\in\mathcal P^{n,p},\theta \geq 0}\quad &\theta
 \tag{LR-D-LB$_r$}\label{mod:cutting-plane-master}\\
 \textnormal{s.t.}\quad
 &\theta\geq h(\mathbf P_\ell)
  \bigl(1-\Delta_{\mathbf P_\ell}(\mathbf P)\bigr),\quad \forall\ell\in[r].\nonumber
\end{align}
Notice that, just as with $g(\mathbf P)$, we can exploit the target-wise separability of $h(\mathbf P)$ to generate cuts, as described in \eqref{eq:sep-cuts}.

The same target-disaggregated master applies to \eqref{mod:mpco-aggregated}, because
its block-diagonal construction preserves the separability of the fixed-projector subproblem.  Its block-diagonal PSD constraint is equivalently
represented by a collection of target-wise PSD blocks.
The master may be solved with the rank-one projector factorization, as in
\cite{BertsimasCoryWrightPauphilet2022}.  Alternatively, a master solver
may use the equivalent incidence representation $\z_j=P_j\e_n$ for the finite
projector set and reconstruct
$P_j=\z_j\z_j^T/(\e_n^T\z_j)$ at every incumbent.  A conic solver then
evaluates $h(\mathbf P)$, and
\eqref{ctr:optimality-cut} is returned to the master.  The incidence vectors
are only an exact encoding of the partition projectors in the master.  They
are not coupled to distances and introduce no big-$M$ inequalities.
\begin{algorithm}[H]
\caption{Outer Approximation for \eqref{mod:mpco-disaggregated}}
\label{alg:cutting-plane}
\begin{algorithmic}[1]
\Require Initial $\mathbf P_1\in\mathcal P^{n,p}$, tolerance $\varepsilon>0$.
\State Set $r\gets1$, $\mathrm{LB}\gets-\infty$, and
$\mathrm{UB}\gets+\infty$.
\Repeat
  \State Solve $h(\mathbf P_r)$ with a conic solver.
  \State Add \eqref{ctr:optimality-cut} and set
  $\mathrm{UB}\gets\min\{\mathrm{UB},h(\mathbf P_r)\}$.
  \State Solve the master problem \eqref{mod:cutting-plane-master} globally;
  let $(\mathbf P_{r+1},\theta_{r+1})$ be an optimal solution.
  \State Set $\mathrm{LB}\gets\theta_{r+1}$ and $r\gets r+1$.
\Until{$\mathrm{UB}-\mathrm{LB}\leq\varepsilon$.}
\State \Return The best feasible projector and its associated clustering solution.
\end{algorithmic}
\end{algorithm}
Every master value is a valid lower bound because all cuts underestimate
$h$, whereas every fixed-projector value is the objective of a feasible
clustering and therefore an upper bound.  Algorithm
\ref{alg:cutting-plane} consequently maintains a valid optimality gap.
\begin{theorem}
If every fixed-projector subproblem is solved exactly and the master problem
is solved globally at every iteration, then Algorithm
\ref{alg:cutting-plane} terminates after finitely many iterations with a
globally optimal solution of \eqref{mod:mpco-disaggregated}.
\end{theorem}

\begin{proof}
Cut \eqref{ctr:optimality-cut} reproduces the exact value at its evaluated
projector and has a nonpositive right-hand side at every other labeled
partition.  Hence, every master value is a valid lower bound and every
fixed-projector value is an upper bound.  The set $\mathcal P^{n,p}$ is
finite.  Each iteration either closes the optimality gap or evaluates a new
projector.  Finite termination follows.
\end{proof}

For formulation \eqref{mod:mpco-sum-aggregated}, let
$
 \PP_{\rm SA}^{n,p}:=
 \left\{P\in\PP_p^n\cap\R_+^{n\times n}:P\e_n=\e_n\right\}.
$
This finite set represents the feasible partitions without distinguishing
their target labels.  For a fixed $P\in\PP_{\rm SA}^{n,p}$, define
$$
 \widetilde{h}(P):=
 \min_{Q,\bm t}\left\{
 \sum_{i=1}^n\omega_i t_i:
 t_i\geq d(\ab_i,Q^T\bm\epsilon_i),\; i\in[n],\; Q=PQ
 \right\}.
$$

For $P\in\PP_{\rm SA}^{n,p}$, write $c_{ih}(P)=1$ if $P_{ih}>0$, and zero
otherwise.  Thus, $c_{ih}(P)$ indicates whether $i$ and $h$ are in the same
cluster.  For an evaluated projector $\bar P$, define
$$
 \Delta_{\bar P}(P):=
 \sum_{\substack{i<h:\,c_{ih}(\bar P)=1}}\bigl(1-c_{ih}(P)\bigr)
 +\sum_{\substack{i<h:\,c_{ih}(\bar P)=0}}c_{ih}(P).
$$
This label-invariant quantity is zero if $P=\bar P$ and at least one otherwise.
Since $\widetilde h(\bar P)\geq0$, the corresponding logic cut is
$
 \theta\geq
 \widetilde h(\bar P)\bigl(1-\Delta_{\bar P}(P)\bigr).
$
It is exact at $\bar P$ and otherwise is dominated by $\theta\geq0$.  In the
discrete master, $c_{ih}(P)$ is a binary coassignment variable.  After evaluating
$P_1,\ldots,P_r$, the master reads
\begin{align}
 \min_{P\in\PP_{\rm SA}^{n,p},\theta\geq 0}\quad &\theta
 \tag{LR-SA-LB$_r$}\label{mod:mpco-sa-cutting-plane-master-short}\\
 \textnormal{s.t.}\quad
 &\theta\geq \widetilde{h}(P_\ell)
 \bigl(1-\Delta_{P_\ell}(P)\bigr),\quad \ell\in[r].\nonumber
\end{align}
Note that if the fixed-projector subproblems and the master problems are solved exactly,
Algorithm~\ref{alg:cutting-plane}, with $\mathbf P$ and $h$ replaced by $P$
and $\widetilde{h}$, terminates finitely with a globally optimal solution of
\eqref{mod:mpco-sum-aggregated}.

Indeed, each iteration either closes the optimality gap or evaluates a new
member of the finite set $\PP_{\rm SA}^{n,p}$.  This variant uses one
rank-$p$ projector and removes the symmetry associated with target labels.

\section{Strengthening and Implementation Strategies}
\label{sec:strengthening}
The exact decomposition methods can be strengthened at two different
levels.  Valid inequalities and semidefinite bounds improve the lower
approximation represented by the discrete master.  Warm starts and the
management of the master--subproblem interaction improve the practical execution
without changing the exactness of the methods.  We present these two levels
separately.

\subsection{Valid Inequalities and Relaxation Strengthening}
We first collect the inequalities that are valid for every feasible
projector.  They can be imposed independently of the initialization and of
the sequence of partitions generated by the algorithms, and they strengthen
both the root relaxation and the subsequent master problems.

\paragraph{Semidefinite Projection Bound.}
Let $
 D^\omega_{ik}:=\min\{\omega_i,\omega_k\}d(\ab_i,\ab_k),$
 for $i,k\in[n],
$ define a weighted distance matrix. The triangle inequality implies
$D^\omega_{ik}\leq\omega_i d(\ab_i,\x)+
\omega_kd(\ab_k,\x)$ for every $\x$.  Replacing $Y(\x_j)$ by $D^\omega$
and relaxing idempotency gives the disaggregated relaxation
\begin{align}
 \min_{P_1,\ldots,P_p}\quad
 &\frac12\sum_{j=1}^p\langle D^\omega,P_j\rangle
 \tag{SDP-RELAX-D}\label{mod:disaggregated-sdp-bound}\\
 \textnormal{s.t.}\quad
 &\sum_{j=1}^pP_j\e_n=\e_n,\nonumber\\
 &\tr(P_j)=1, &&\forall j\in[p],\nonumber\\
 &P_j\succeq0,\quad P_j\geq0, &&\forall j\in[p].\nonumber
\end{align}
Every feasible solution of \eqref{mod:disaggregated-sdp-bound} yields an
aggregate matrix $P=\sum_jP_j$.  The corresponding aggregate problem has a single matrix variable.
\begin{align}
 \min_{P\succeq0, P\geq0}\quad
 &\frac12\langle D^\omega,P\rangle
 \tag{SDP-RELAX}\label{mod:scalable-sdp-bound}\\
 \textnormal{s.t.}\quad
 &P\e_n=\e_n,\nonumber\\
 &\tr(P)=p.\nonumber
\end{align}
The two basic relaxation bounds are equal. For any feasible aggregate
matrix $P$, the split $P_j=P/p$ is feasible in
\eqref{mod:disaggregated-sdp-bound} and preserves the objective.
Additional target-specific constraints can break this equivalence.  The disaggregated model
also provides target-wise fractional assignment scores for rounding. 
Semidefinite relaxations and rounding are well-established tools for
partition problems
\citep{Goemans1997semidefinite,LuoMaSoYeZhang2010,peng2007approximating}.
\paragraph{Projection Inequalities.}
For the aggregate partition projector $P=\sum_{j=1}^pP_j$, one has
$P_{ih}=|C|^{-1}$ whenever points $i$ and $h$ belong to the same cluster $C$, and
$P_{ih}=0$ otherwise.  In particular, $P_{ih}\leq P_{ii}$ states that the
normalized association between two points cannot exceed the diagonal value
of the cluster containing $i$.  The inequalities
$
 P_{ih}+P_{ik}\leq P_{ii}+P_{hk},$ for $i,h,k\in[n],
$
provide a normalized transitivity condition.  For example, if $i$ is
associated with both $h$ and $k$, then $h$ and $k$ must also be associated.
Both families are satisfied by every partition projector and exclude
fractional eigenvalues admitted by the basic doubly nonnegative relaxation.
They cannot be imposed directly on the \ref{eq:weber_original} formulation using only the
$z$-variables, because
$P_{ih}=\sum_{j=1}^p z_{ij}z_{hj}/\bigl(\sum_{\ell=1}^n z_{\ell j}\bigr)$ involves binary products and
the cluster cardinalities.  Their direct linear form is therefore a specific
advantage of working in the projector space.  An equivalent extension of
\ref{eq:weber_original} would require additional coassignment and cardinality variables.

\paragraph{Perspective and Trace Bounds.} The mixed-projection relaxations can combine distance epigraphs with
$$
 \begin{pmatrix}P&Q\\Q^T&M\end{pmatrix}\succeq0,
$$
or the corresponding target-wise constraints.  A valid bound derived
from a bounded region containing optimal centers is essential.  Let
$\bar{\ab}:=\frac{1}{n}\sum_{i=1}^n\ab_i$ and translate the origins by
$\bar{\ab}$.  Translation does not change the distances or the optimal
value.  Let $\mathcal X\subseteq\R^m$ be any bounded region known to
contain an optimal center for every possible nonempty cluster.  It is enough
to choose
\begin{equation}\label{eq:tracebound}
     R^2\geq\sup\{\langle\x-\bar{\ab},\x-\bar{\ab}\rangle :\x\in\mathcal X\}.
\end{equation}
Thus, the validity of the trace bound does not depend on a particular norm or
on a particular construction of $\mathcal X$.  At an exact solution, the trace represents the sum of the squared norms of the
centers, repeated according to their cluster cardinalities.  Since these
cardinalities sum to $n$, the following bounds are valid:
$
\tr(M)\le n R^2$ for \eqref{mod:mpco-sum-aggregated} and $
\sum_{j=1}^p\tr(M_j)\le nR^2$
 for \eqref{mod:mpco-disaggregated}/\eqref{mod:mpco-aggregated}.
These bounds remain valid for the exact formulations and restrict the
otherwise weak perspective relaxations.

\paragraph{Certified SDP Root Bound.}
Either SDP relaxation can also be used independently of rounding.  Let
$L_{\rm SDP}$ denote a certified lower bound obtained from the selected model.  When
the SDP is stopped before optimality, $L_{\rm SDP}$ must be taken from a
dual-feasible solution rather than from the available primal objective.  The lower bound requirement $\theta\geq L_{\rm SDP}$ 
can be imposed in \eqref{mod:pairwise-benders-master} and \eqref{mod:mpco-sa-cutting-plane-master-short}.  For \eqref{mod:cutting-plane-master}, the corresponding
inequality is $\sum_{j=1}^p\theta_j\geq L_{\rm SDP}$.  Hence, the same SDP
calculation produces an initial partition and strengthens the root lower
bound.

\paragraph{Projector Distance Bounds.}
The distance lower bound can be imposed directly in the projection space,
without introducing assignment or coassignment variables.  For every
feasible cluster projector $P_j$ and every target location $\x_j$, the
triangle inequality gives
$
 \langle D^\omega,P_j\rangle
 \leq \langle Y(\x_j),P_j\rangle.
$
Consequently, the pairwise master \eqref{mod:pairwise-benders-master} can be strengthened with
$
 \theta\geq\frac12\sum_{j=1}^p\langle D^\omega,P_j\rangle.
$
This inequality is linear in the projector entries.  Minimizing it over the
doubly nonnegative relaxation of $\mathcal P^{n,p}$ gives
\eqref{mod:scalable-sdp-bound}.  Thus, this SDP is the natural relaxation of
the projector bound above, rather than a
bound imported from \eqref{eq:weber_original}.

The same argument applies to the low-rank models. The disaggregated master \eqref{mod:cutting-plane-master} with epigraph variables $\theta_j$ can use
the projection inequalities
$
 \theta_j\geq\frac12\langle D^\omega,P_j\rangle,$
 for $j\in[p].
$
For the aggregated formulation, the right-hand sides are summed over the
diagonal projector blocks.  These bounds complement the dual supporting
cuts \eqref{ctr:optimality-cut}.  The former use the geometry of pairwise
distances, while the latter use the complete fixed-projector conic
subproblem.  Both families are expressed solely through the projection
variables that define the formulations.

\paragraph{Symmetry Breaking.}
Finally, the interchangeability of target labels creates equivalent
partitions.  One may use the following inexpensive conditions to partially
break this symmetry.  If origin $i$ is not assigned to one of facilities
$1,\ldots,i-1$, then it is assigned to target $i$.  In terms of the
assignment variables, this implication is
\begin{equation*}
 z_{ii}\geq1-\sum_{j=1}^{i-1}z_{ij},\qquad i\in[p].
\end{equation*}
Since every origin is assigned to exactly one target, the above set of constraints is equivalently written as
$\sum_{j=1}^i z_{ij}=1$ for $i\in[p]$.  In particular, it implies
$z_{11}=1$ and prevents origin $i$ from using a target with label
larger than $i$.

For the projector-based formulations, $z_{ij}=\bm\epsilon_i^TP_j\e_n$.  Therefore, the same
conditions can be imposed directly on the projectors as
$
 \sum_{j=1}^i \bm \epsilon_i^TP_j\e_n=1,$ for $i\in[p].
$
Every unlabeled nonempty partition has a labeling satisfying these
conditions, for example by ordering its clusters according to the first
origin appearing in each cluster.  Consequently, the conditions do
not change the optimal value.  They remove part, but not necessarily all, of
the target-label symmetry.  Their main advantage is that they add only
$p$ linear constraints and require no auxiliary variables.  A warm-start
partition can be relabeled by first appearance before it is passed to the
master. A stronger option is the complete first-occurrence ordering
$$
z_{ij}\leq\sum_{h=1}^{i-1}z_{h,j-1},\qquad i\in[n],\quad j=2,\ldots,p,
$$
or, equivalently,
$\bm\epsilon_i^TP_j\e_n\leq\sum_{h=1}^{i-1}
\bm\epsilon_h^TP_{j-1}\e_n$.  This selects the representative in which
clusters are ordered by their first origin and removes the remaining
target-label symmetry from the labeled masters.

\subsection{Warm Starts and Algorithmic Implementation}
The preceding inequalities improve lower bounds.  We now describe how
feasible solutions are initialized and how the exact master--subproblem
decompositions are executed.  These choices affect computational performance
but do not replace the lower-bound certificates required for exactness.

The exact methods developed in Sections \ref{sec:pairwise} and
\ref{sec:mpco} can start from any feasible nonempty partition.  Their
performance, however, depends strongly on the quality of the initial
incumbent.  A good partition provides an early upper bound, guides the search
toward promising projector configurations, and may reduce the number of
conic subproblems evaluated before the optimality gap closes.  Warm starts
are particularly useful because the assignment and location decisions
reinforce each other.

This observation is the basis of classical location-allocation heuristics
and their modern variants
\citep{BongartzCalamaiConn1994,BrimbergHansenMladenovicTaillard2000,Cooper1964}.
It also motivates the use of several starting configurations rather than a
single deterministic initialization \citep{DreznerBrimberg2026}.  Warm starts
produce feasible solutions and upper bounds only.  Global optimality is still
certified by the exact methods.

\paragraph{Multiple Starting Partitions.}
Candidate target locations can be selected from the origins or from
geometrically dispersed subsets of them.  Each origin is assigned to a
closest candidate target, with empty facilities repaired by moving origins from non-singleton clusters.  The resulting partitions provide a
diverse collection of inexpensive initial solutions.  This option follows
the general multistart principle used in distance-based assignment.
\citep{Cooper1964,DreznerBrimberg2026}.

\paragraph{Alternating Location and Allocation.}
Starting from a feasible partition, fix its projectors.  Problem
\eqref{mod:pairwise} then separates into the $p$ single-target problems
$$
 \min_{\x_j\in\R^m}\ \frac12\langle P_j,Y(\x_j)\rangle,
 \qquad j\in[p].
$$
After computing the target locations, each origin is reassigned to a
closest target while preserving nonempty clusters.  These two steps are
repeated until the assignment does not change or the objective improvement
is below a tolerance.  The objective is nonincreasing and every iterate is
feasible.  This is the projector version of the classical location-allocation
scheme \citep{BongartzCalamaiConn1994,Cooper1963,Cooper1964}.   The
single-target problems can be handled by any suitable conic solver.  The
allocation update assigns each point directly to a closest target, with a
repair step if a cluster becomes empty.
\paragraph{SDP Rounding.}
Let $(P_1^\star,\ldots,P_p^\star)$ solve
\eqref{mod:disaggregated-sdp-bound}, and define the fractional score matrix
$\Lambda\in\R^{n\times p}$ by
$\lambda_{ij}:=\bm\epsilon_i^TP_j^\star\e_n$.  A nonempty partition is
recovered from
\begin{equation*} 
 \max_{Z\in\{0,1\}^{n\times p}}\{\langle\Lambda, Z\rangle :
 Z\e_p = \e_n, \; Z^T\e_n \ge \e_p\}.
\end{equation*}
The rounded partition is refined by alternating location and allocation.  Its
objective value initializes the upper bound, its projectors initialize the
exact master problem, and the SDP value remains available as an independent
lower bound.
\paragraph{Certified Subproblem Values and Cut Management.}
For a partition proposed by the discrete master, a conic solver handles the
associated convex location subproblems.  A feasible center and its primal cost give a
valid solution and can update the incumbent.  In contrast, only a
dual-feasible objective is used as a certified lower value.  If $L(C)$ is the
dual bound obtained for a cluster $C$, one may cache this value and add,
for every target label $j$, the valid cut
$$
 \theta_j\geq L(C)\left(\sum_{i\in C}z_{ij}-|C|+1\right).
$$
The cut is active when target $j$ contains all points of $C$.  The certified
values of all clusters in the current partition are also combined in a
logic-based optimality cut that excludes repeating that projector without
accounting for its certified cost.  Such a cache permits both cut types to
be reused for other target labels and later partitions.  If the conic solver
returns a primal feasible center but no dual-feasible bound, the primal value may still
improve the incumbent, but it is never used on the right-hand side of a
lower-bounding cut.  If a certificate required by the exact step is missing,
the evaluation remains unresolved and the algorithm reports a subproblem
limit instead of claiming global optimality.

The strategies described above can be combined in different ways.
Section~\ref{sec:experiments} specifies the configuration adopted for each
formulation in our computational study.

\section{Numerical Experiments}
\label{sec:experiments}

This section compares the exact formulations and studies whether the
projection models improve the lower bounds and the final optimality gaps.
We report solution time together with the percentage of unsolved instances,
since a smaller time caused by early termination must not be interpreted as
better performance.

\subsection{Instances and Setting}
We use the planar Eilon instances~\citep{EilonWatsonGandyChristofides1971},
available in the repository: \cite{Blanco2020MCLPIF}.
The test set contains five instances for each
$n\in\{10,20,30,40\}$, as well as the complete instance for $n=50$. We consider $p\in\{2,3,5\}$ facilities in three different $\ell_\tau$ metric spaces, with $\tau\in\{3/2,2,3\}$.  We compare the
big-$M$ integer conic formulation, denoted by \ref{eq:weber_original}, the pairwise model \ref{mod:pairwise}, and the
disaggregated, aggregated, and sum-aggregated low-rank models \ref{mod:mpco-disaggregated}, \ref{mod:mpco-aggregated}, and
\ref{mod:mpco-sum-aggregated}.  Gurobi 13.0.1 solves
\ref{eq:weber_original} and the discrete masters, while MOSEK 11.0.29 solves
the SDP relaxations and conic subproblems.  The experiments were run on a
2022 Mac Studio with an Apple M1 Max processor, 64 GB of memory, and macOS
Tahoe 26.6.2.  The time
limit was 3600 seconds, the root SDP was
limited to 300 seconds, and the relative optimality tolerance was $10^{-4}$. For every instance, the common warm start is computed once before the
five methods are called, so its construction time is not included in their
individual runtimes.  The reported time of each projection method includes
its method-specific root SDP, the strengthened Gurobi master, and all
subsequent MOSEK cluster evaluations. 
The code, processed instances, complete run-level results, and scripts used to
generate the tables and figures are available in the associated repository:
\cite{BlancoLabordaGutierrezMartinezAnton2026Projectors}.

All methods use the best warm start from ten initial partitions, alternating
location-allocation, and SDP rounding, together with the SDP lower bound and
prefix symmetry.  The projection models also use complete ordering,
projection inequalities, geometric and subset cuts, and dual-certified cuts.
\ref{mod:pairwise} uses the pairwise relaxation.  \ref{mod:mpco-disaggregated} and \ref{mod:mpco-aggregated} use target-wise
perspectives and the aggregate trace bound, while \ref{mod:mpco-sum-aggregated} uses one aggregated
perspective and trace bound.  These are the complete strategies of
Section~\ref{sec:strengthening}.
For the trace bounds: if $\tau=2$ we choose
$R^2$ as \eqref{eq:tracebound} for $\mathcal X=\{\ab_i\}_{i\in[n]}$.  In this case, the
convex hull of the origins is contained in the corresponding Euclidean
ball.  For $\tau\in\{3/2,3\}$, we use
$R^2=\sum_{r=1}^m(\max_{1\leq i\leq n}|a_{ir}-\bar a_r|)^2$.  This is the squared radius
of the coordinate-wise box containing the centered origins.  These
regions contain an optimal center for every cluster and therefore satisfy the
general validity condition stated in Section~\ref{sec:strengthening}.
In the Gurobi masters, the labeled methods use one epigraph variable per
target and the reusable subset cuts described above.  Formulation \ref{mod:mpco-sum-aggregated} keeps a
canonical incidence encoding only to construct its master constraints, while
its partition logic cut is written in the label-invariant coassignment
variables.  The projector methods use the incidence representation in the
master and reconstruct the corresponding projector at each incumbent.  We
use an explicit solve-separate-resolve loop: Gurobi solves the master,
MOSEK certifies the fixed-partition conic subproblems, and the resulting cuts
are added before the master is solved again.

The numerical study reports both primal and dual effects of these
strategies.  For each formulation, we record the best feasible value, the
initial continuous-relaxation bound, the SDP bound, the final certified
lower bound, and the computing time required by each component.  Relaxation
quality is measured by
$
 \frac{U-L}{\max\{1,|U|\}},
$
where $U$ is the best feasible value and $L$ is the bound under comparison.
For every method, $L$ is the continuous relaxation of the strengthened
implemented master after adding its symmetry constraints and the available
certified SDP root bound.  For the projector methods, the method-specific SDP
bound is also stored separately.  Thus, the reported continuous-relaxation
gap compares the implemented strengthened masters, rather than the raw
continuous relaxations of the algebraic formulations.  This distinction
allows us to assess whether an improvement comes from a better feasible
solution or from a stronger lower bound.

\subsection{Computational Results} Table~\ref{tab:eilon-main-results} reports averages over the common runs
available for all five methods.  Runtime, final optimality gap, and percentage of
unsolved instances are displayed in separate blocks.  The best value of each
measure for every $(n,\tau)$ is in bold.  A run is classified as unsolved
when the corresponding method does not certify optimality within the
one-hour limit.  Only runs available for every method are included, so all
entries within an $(n,\tau)$ row are based on the same instances.

\begin{table}[htbp]
\caption{Average runtime, final optimality gap, and percentage of unsolved instances. Boldfaced values indicate the best performance within each $(n, \tau)$ group (except those reaching timelimit for all instances or 100\% final GAP).}
\label{tab:eilon-main-results}
\centering
\tiny
\setlength{\tabcolsep}{2.2pt}
\resizebox{\textwidth}{!}{%
\begin{tabular}{cc||rrrrr||rrrrr||rrrrr}
\hline\hline
& & \multicolumn{5}{c||}{Time (s)} & \multicolumn{5}{c||}{Opt. Gap (\%)} & \multicolumn{5}{c}{Unsolved (\%)} \\
$n$ & $\tau$ & \ref{eq:weber_original} & \ref{mod:pairwise} & \ref{mod:mpco-disaggregated} & \ref{mod:mpco-aggregated} & \ref{mod:mpco-sum-aggregated} & \ref{eq:weber_original} & \ref{mod:pairwise} & \ref{mod:mpco-disaggregated} & \ref{mod:mpco-aggregated} & \ref{mod:mpco-sum-aggregated} & \ref{eq:weber_original} & \ref{mod:pairwise} & \ref{mod:mpco-disaggregated} & \ref{mod:mpco-aggregated} & \ref{mod:mpco-sum-aggregated} \\
\hline\hline
\multirow{3}{*}{$10$} & $3/2$ & 1.1 & \textbf{0.9} & 2.0 & 1.9 & 1.4 & 0.0 & 0.0 & \textbf{0.0} & 0.0 & 0.0 & \textbf{0.0} & \textbf{0.0} & \textbf{0.0} & \textbf{0.0} & \textbf{0.0} \\
 & $2$ & \textbf{0.4} & 0.6 & 1.2 & 1.1 & 1.0 & 0.0 & 0.0 & \textbf{0.0} & 0.0 & 0.0 & \textbf{0.0} & \textbf{0.0} & \textbf{0.0} & \textbf{0.0} & \textbf{0.0} \\
 & $3$ & 1.2 & \textbf{0.8} & 2.0 & 1.8 & 1.4 & 0.0 & 0.0 & 0.0 & \textbf{0.0} & 0.0 & \textbf{0.0} & \textbf{0.0} & \textbf{0.0} & \textbf{0.0} & \textbf{0.0} \\
\multirow{3}{*}{$20$} & $3/2$ & 1454.2 & \textbf{200.4} & 258.0 & 284.2 & 1094.5 & 8.5 & \textbf{0.0} & 0.0 & 0.0 & 0.2 & 40.0 & \textbf{0.0} & \textbf{0.0} & \textbf{0.0} & 13.3 \\
 & $2$ & 869.2 & \textbf{241.5} & 270.0 & 263.7 & 1201.5 & 2.5 & \textbf{0.0} & \textbf{0.0} & \textbf{0.0} & 0.1 & 13.3 & \textbf{0.0} & \textbf{0.0} & \textbf{0.0} & 6.7 \\
 & $3$ & 1436.7 & \textbf{242.7} & 269.8 & 245.5 & 1191.3 & 9.7 & \textbf{0.0} & 0.0 & 0.0 & 0.2 & 33.3 & \textbf{0.0} & \textbf{0.0} & \textbf{0.0} & 20.0 \\
\multirow{3}{*}{$30$} & $3/2$ & 2905.9 & \textbf{2745.2} & 2778.0 & 2775.1 & 3117.6 & 24.1 & \textbf{17.0} & 18.7 & 17.2 & 22.4 & \textbf{66.7} & 73.3 & 73.3 & 73.3 & 73.3 \\
 & $2$ & \textbf{2581.4} & 2796.0 & 2830.7 & 2818.2 & 3228.0 & 23.5 & 19.0 & \textbf{18.9} & 20.4 & 22.5 & \textbf{66.7} & 73.3 & 73.3 & 73.3 & 86.7 \\
 & $3$ & 2941.9 & \textbf{2818.1} & 2830.9 & 2846.7 & 3286.2 & 23.9 & 17.9 & 17.4 & \textbf{17.3} & 22.1 & \textbf{66.7} & 73.3 & 73.3 & 73.3 & 80.0 \\
\multirow{3}{*}{$40$} & $3/2$ & 3600.0 & 3558.0 & \textbf{3553.1} & 3598.5 & 3600.0 & 34.8 & 23.1 & 23.0 & \textbf{22.9} & 26.6 & 100.0 & \textbf{93.3} & \textbf{93.3} & \textbf{93.3} & 100.0 \\
 & $2$ & 3600.0 & 3600.0 & 3600.0 & 3600.0 & 3600.0 & 33.3 & \textbf{24.8} & 25.0 & 25.0 & 27.9 & 100.0 & 100.0 & 100.0 & 100.0 & 100.0 \\
 & $3$ & 3600.0 & 3600.0 & 3600.0 & 3600.0 & 3600.0 & 34.8 & \textbf{23.6} & 23.9 & 23.7 & 27.1 & 100.0 & 100.0 & 100.0 & 100.0 & 100.0 \\
\multirow{3}{*}{$50$} & $3/2$ & 3600.0 & 3600.0 & 3600.0 & 3600.0 & 3600.0 & 35.7 & \textbf{28.0} & 31.1 & 29.8 & 32.5 & 100.0 & 100.0 & 100.0 & 100.0 & 100.0 \\
 & $2$ & 3600.0 & 3600.0 & 3600.0 & 3600.0 & 3600.0 & 35.6 & \textbf{30.7} & 35.5 & 33.7 & 34.9 & 100.0 & 100.0 & 100.0 & 100.0 & 100.0 \\
 & $3$ & 3600.0 & 3600.0 & 3600.0 & 3600.0 & 3600.0 & 35.4 & \textbf{28.7} & 30.5 & 30.6 & 34.3 & 100.0 & 100.0 & 100.0 & 100.0 & 100.0 \\
\hline\hline
\end{tabular}
}
\end{table}

All methods solve the $n=10$ instances in about two seconds on average.  For
$n=20$, formulations \ref{mod:pairwise}, \ref{mod:mpco-disaggregated}, and
\ref{mod:mpco-aggregated} solve every run, with \ref{mod:pairwise} requiring only
$200.4$--$242.7$ seconds.  Formulation \ref{eq:weber_original} leaves $13.3\%$--$40.0\%$ unsolved, while the
rate for \ref{mod:mpco-sum-aggregated} is $6.7\%$--$20.0\%$.  At $n=30$, \ref{eq:weber_original} certifies one more run
than \ref{mod:pairwise}, \ref{mod:mpco-disaggregated}, and \ref{mod:mpco-aggregated}, but their average gaps are only $17.0\%$--$20.4\%$,
compared with $23.5\%$--$24.1\%$ for \ref{eq:weber_original}.  For $n=40$, \ref{mod:pairwise}, \ref{mod:mpco-disaggregated}, and \ref{mod:mpco-aggregated}
each certify one run for $\tau=3/2$, and no method certifies a run for
$\tau\in\{2,3\}$.  Formulation \ref{mod:mpco-aggregated} gives the smallest gap for $\tau=3/2$, and \ref{mod:pairwise} for
$\tau\in\{2,3\}$.  At $n=50$, all runs reach the time limit and \ref{mod:pairwise} gives the
smallest gap for every norm.  Hence, the main benefit on the largest
instances is a stronger final certificate rather than a shorter runtime.

\begin{figure}[htbp]
\centering
\includegraphics[width=.72\textwidth]{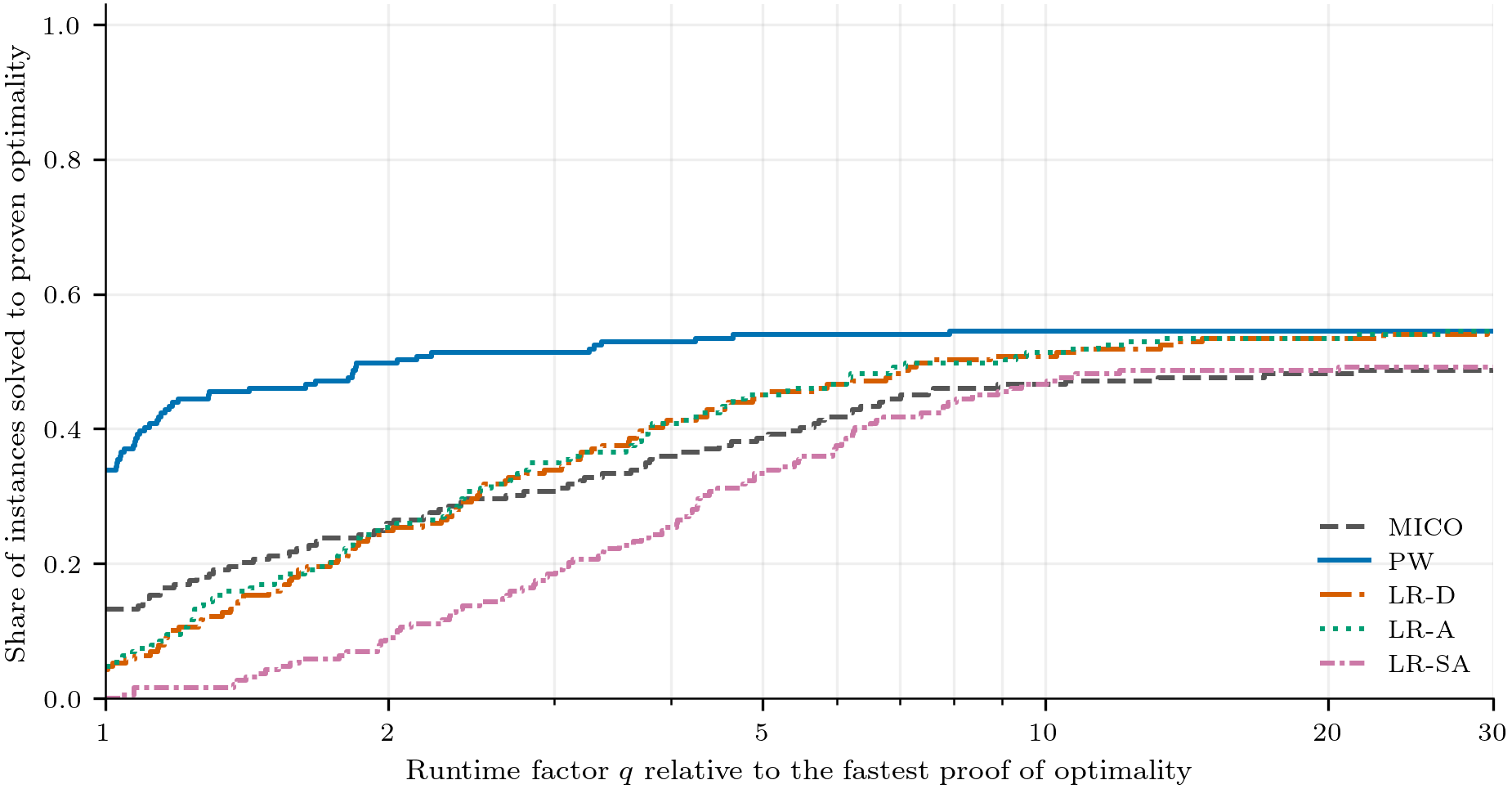}
\caption{Performance profiles over all matched instances and tested sizes.}
\label{fig:pp}
\end{figure}
In Figure~\ref{fig:pp}, the horizontal value $q$ is the runtime relative to
the fastest proof of optimality, and the profile gives the fraction solved
within this factor.  Formulation \ref{mod:pairwise} dominates for small $q$ and reaches its final
optimality rate much earlier.  Formulations \ref{mod:pairwise}, \ref{mod:mpco-disaggregated}, and \ref{mod:mpco-aggregated} attain the same final
optimality rate, but \ref{mod:mpco-disaggregated} and \ref{mod:mpco-aggregated} require larger runtime factors, while \ref{eq:weber_original}
and \ref{mod:mpco-sum-aggregated} prove optimality for a smaller fraction of the instances.  Thus,
\ref{mod:pairwise} gives the best balance between speed and
the ability to prove optimality.
\begin{figure}[htbp]
\centering
\includegraphics[width=0.7\textwidth]{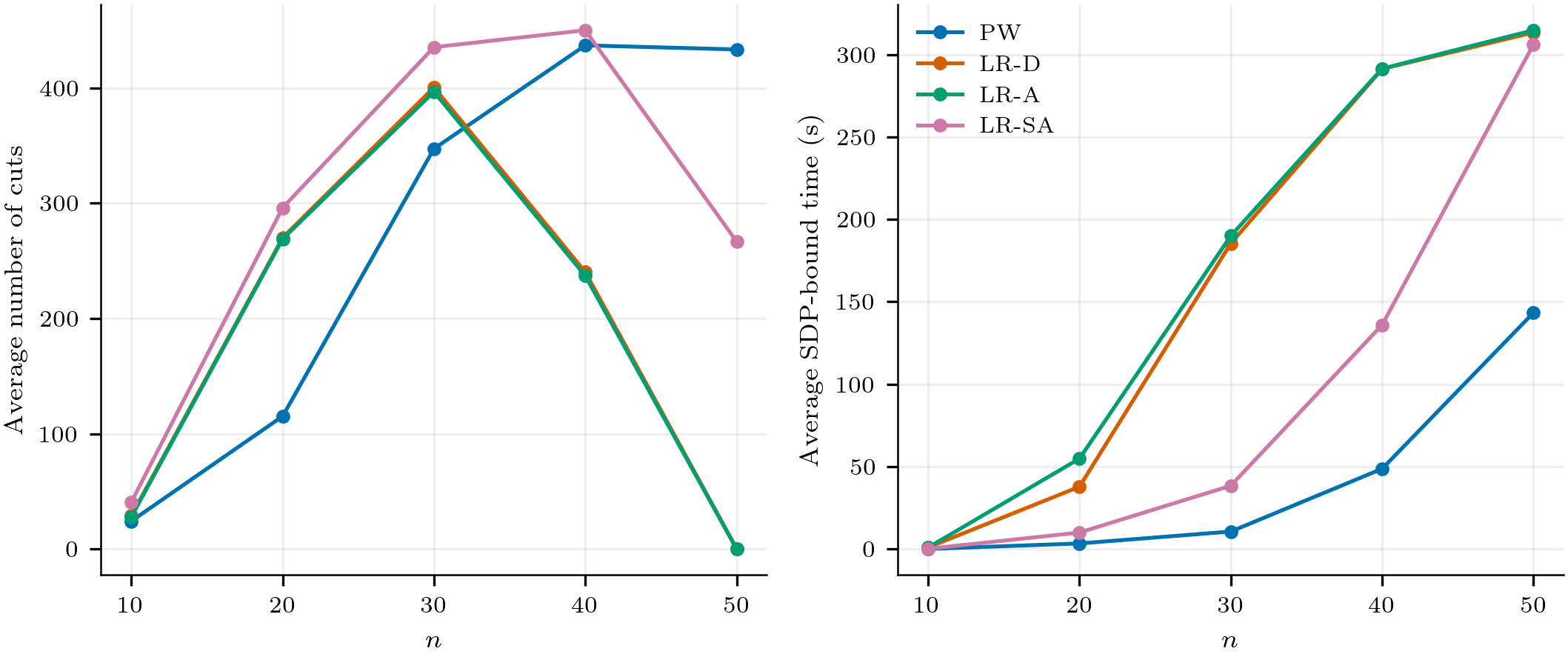}
\caption{Average number of partition-logic and separated SDP triangle
cuts generated by the projection formulations (left) and average time spent computing SDP bounds, classified by
instance size (right).}
\label{fig:projector-cuts-sdp-time}
\end{figure}

Figure~\ref{fig:projector-cuts-sdp-time} helps explain these computational
differences.  The left panel counts the two cut families recorded in the
experiment files. Cached subset cuts and the static geometric master
constraints are not included.  Formulations \ref{mod:mpco-disaggregated} and \ref{mod:mpco-aggregated} generate close to $400$ recorded
cuts at $n=30$.  Their sharp decline at $n=50$ reflects the higher cost of
processing the master nodes, which leaves time for fewer master solutions and
cuts, not a reduced need for cuts.  Formulations \ref{mod:pairwise} and \ref{mod:mpco-sum-aggregated} follow different cut-count
patterns.  Moreover, \ref{mod:pairwise} spends the least
time computing SDP bounds.  It requires less than $50$ seconds through $n=40$
and about $140$ seconds for $n=50$, against approximately $300$ seconds for
the target-wise models when $n\geq40$.
Thus, the pairwise formulation combines effective cuts with cheaper SDP
computations, consistently with its runtime and final-gap results.

Figure~\ref{fig:mip-gap-improvement} focuses on the instances for which
\ref{eq:weber_original} does not prove optimality.  Let
$\overline G_{0}$ and $\overline G_m$ denote the average final optimality gaps,
in percent, of \ref{eq:weber_original} and projection method $m$.  The bars
report $\overline G_{0}-\overline G_m$ in percentage points, so positive
values favor the projection method.  All projection formulations reduce the
remaining gap.  The reductions are similar for $n=20$, while the pairwise,
disaggregated, and aggregated models are generally strongest for $n=30$ and
$n=40$.  For $n=50$, \ref{mod:pairwise} gives the largest improvement for all
three values of $\tau$.
\begin{figure}[htbp]
\centering
\includegraphics[width=.8\textwidth]{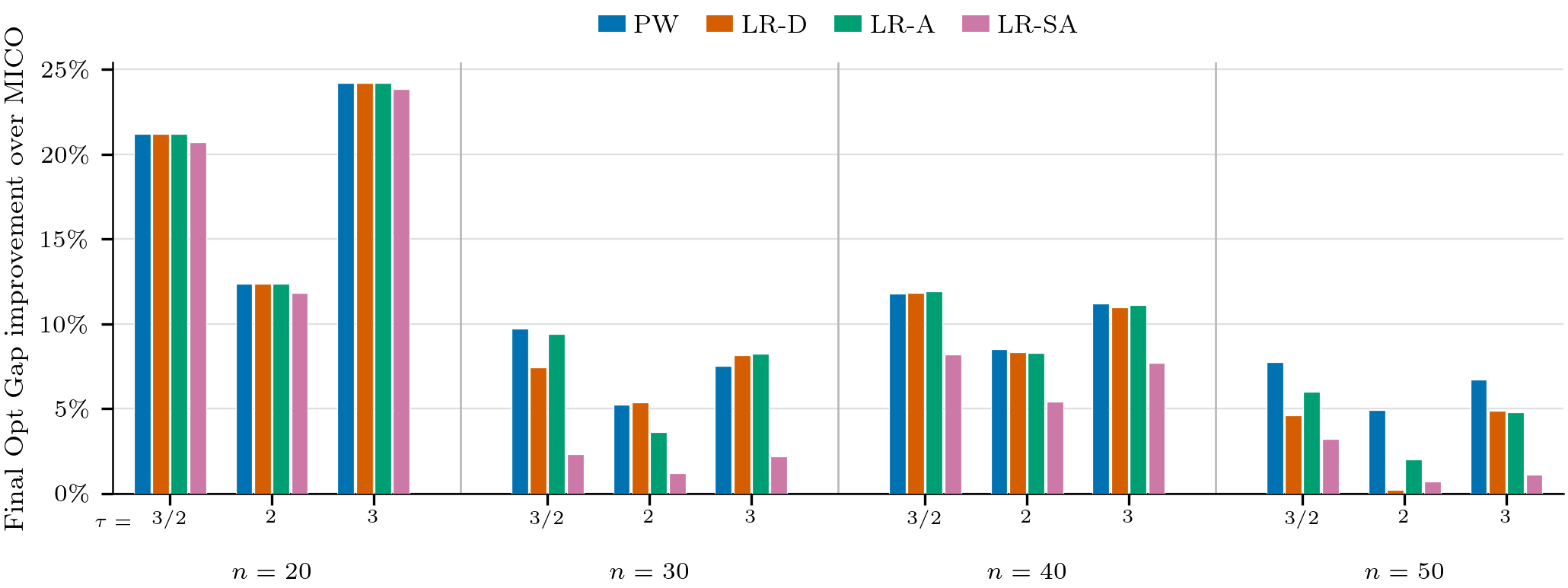}
\caption{Absolute improvement, in percentage points, in the average final optimality gap relative to
\ref{eq:weber_original} on matched instances not solved to proven optimality by
\ref{eq:weber_original}.}
\label{fig:mip-gap-improvement}
\end{figure}

Figure~\ref{fig:relaxation-gap-deviation} compares the continuous-relaxation
LB gap with \ref{eq:weber_original}.  If $g_0$ and $g_m$ denote their
average gaps, the figure reports
$100(g_0-g_m)/g_0$, so positive values indicate a stronger
projection relaxation.
\begin{figure}[htbp]
\centering
\includegraphics[width=.7\textwidth]{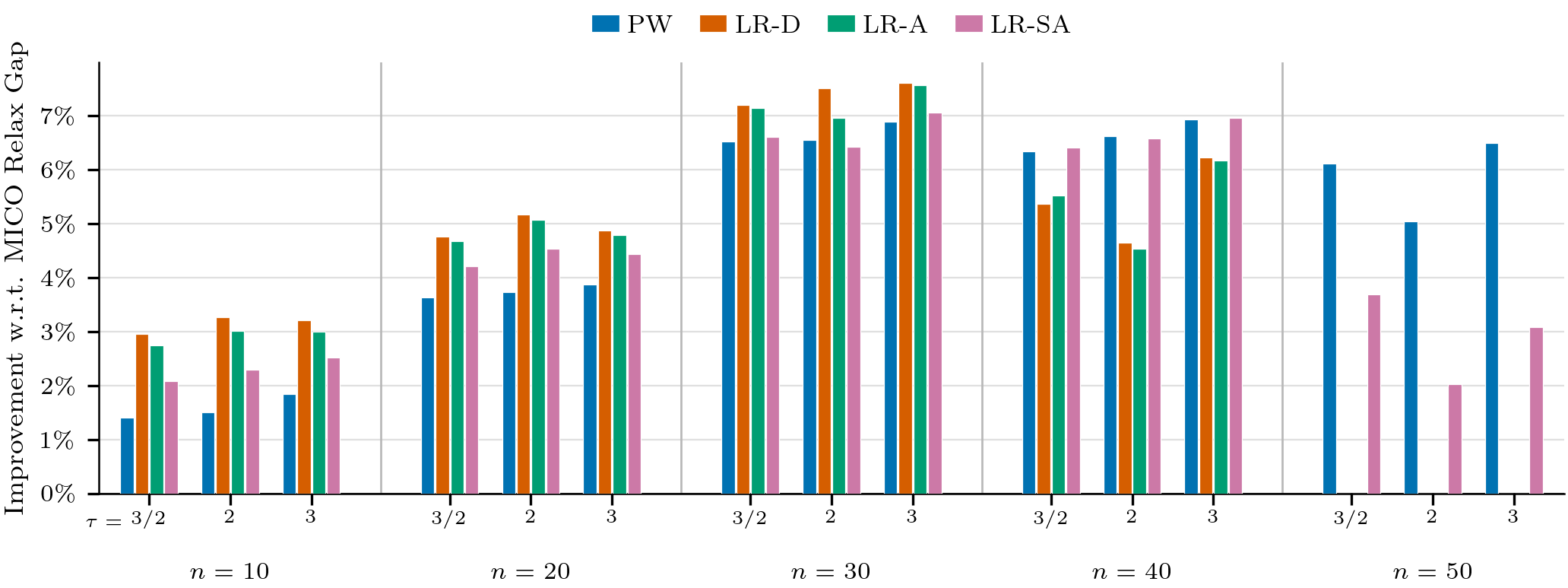}
\caption{Relative improvement of the continuous-relaxation LB gap over \ref{eq:weber_original}.}
\label{fig:relaxation-gap-deviation}
\end{figure}
The projection relaxations are never weaker than \ref{eq:weber_original} in
the matched groups.  Formulations \ref{mod:mpco-disaggregated} and \ref{mod:mpco-aggregated} give the largest gains at $n=30$, while \ref{mod:pairwise}
and \ref{mod:mpco-sum-aggregated} are strongest at $n=40$.  At $n=50$, \ref{mod:pairwise} remains the strongest.
Thus, projection information strengthens the initial lower bounds, although
the best formulation depends on the instance size.

\section{Conclusion}
\label{sec:conclusions}

This paper addressed the distance-based assignment problem through orthogonal projectors and conic optimization formulations. The
resulting pairwise and three low-rank formulations avoid big-$M$ constraints,
and their exact decompositions strengthen discrete masters with certified
conic bounds and cuts.  This provides a new framework for understanding the interplay among assignment, location, rank, and semidefinite optimization.

Computationally, the projection formulations generally strengthen the
continuous relaxation of the classical integer programs.  The pairwise model gives
the best overall balance among runtime, final optimality gap, and proven optimality,
while the low-rank models are competitive in several groups.  These results
are promising because SDP technology for exact decomposition is less mature
than MILP and MISOCP technology, leaving room for further gains in SDP
algorithms, warm starts, and solver interfaces.

Future research includes projector formulations for other
assignment-based problems, matrix-completion techniques for cheaper SDP
relaxations, extensions to other cones such as the exponential cone, and
stronger Benders-type decompositions and certified cuts.


\section*{Code and Data Disclosure}\label{sec:Code and Data Disclosure}

The code and data to support the numerical experiments in this paper can be found at \cite{BlancoLabordaGutierrezMartinezAnton2026Projectors}.

\section*{Acknowledgements}

The authors acknowledge financial support from the Agencia Estatal de Investigaci\'on (AEI) of the Ministerio de Ciencia, Innovación y Universidades (MICIU) through grants: PID2024-156594NB-C21, IMAG--Mar\'ia de Maeztu CEX2020-001105-M and IMUS--Mar\'ia de Maeztu CEX2024-001517-M. They also acknowledge funding from the Programa Andalucía FEDER through project C-EXP-139-UGR23.

\bibliographystyle{plainnat}
\bibliography{refs}

\end{document}